\documentclass[preprint,12pt]{elsarticle}

\usepackage{amsmath}
\usepackage{calc}	

\NewDocumentCommand\bdry{}{\partial}
\NewDocumentCommand\inquotes{m}{``#1''}

\DeclareMathOperator{\interior}{interior}
\DeclareMathOperator{\closure}{closure}
\DeclareMathOperator{\Wh}{Wh}
\DeclareMathOperator{\dist}{dist}
\DeclareMathOperator{\homot}{\simeq}
\DeclareMathOperator{\holink}{holink}
\DeclareMathOperator{\incl}{incl}
\DeclareMathOperator{\desusp}{desusp}
\DeclareMathOperator{\sign}{sign}

\DeclareMathOperator{\cyl}{cyl}
\DeclareMathOperator{\cone}{c}
\DeclareMathOperator{\opencyl}{c\mathring{y}l}
\DeclareMathOperator{\opencone}{\mathring{c}}
\DeclareMathOperator{\im}{im}
\NewDocumentCommand\GlueIntrin{m}{\cup_{(#1|}}

\NewDocumentCommand\Integers{}{\mathbb{Z}}
\NewDocumentCommand\Rationals{}{\mathbb{Q}}
\NewDocumentCommand\Reals{}{\mathbb{R}}
\NewDocumentCommand\txt{m}{\text{\normalfont{#1}}}
\NewDocumentCommand\Bordism{}{\Omega^{\txt{MHSS--Witt}}}
\DeclareMathOperator{\bordant}{\sim}
\DeclareMathOperator{\blank}{--}
\NewDocumentCommand\ptSpace{}{\{\txt{pt}\}}
\DeclareMathOperator{\id}{id}

\NewDocumentCommand\onehalf{}{\frac{1}{2}}
\NewDocumentCommand\SmallSpaceBreakable{m}{\hspace{.16667em}}
\NewDocumentCommand\Slash{}{\,/\SmallSpaceBreakable}
\NewDocumentCommand\ParenthesizedHyphenation{m}{\mbox{(#1-)}\allowbreak}
\NewDocumentCommand\BiCollared{}{\ParenthesizedHyphenation{bi}collared}
\NewDocumentCommand\BiCollar{}{\ParenthesizedHyphenation{bi}collar}
\NewDocumentCommand\defName{m}{\emph{#1}}

\usepackage{amssymb}
\usepackage{amsthm}

\journal{Topology and its Applications}

\begin{document}
	
	\begin{frontmatter}
		
		
		
		\title{Topological Bordism of Singular Spaces and an Application to Stratified L-Classes}
		
		
		\author[label1]{Martin Rabel}
		
		\affiliation[label1]{organization={Mathematisches Institut, Universität Heidelberg},
			city={Heidelberg},
			country={Germany}}
		
		\begin{abstract}
			A generalized-homology bordism-theory is constructed, such that for certain manifold homotopy stratified sets (MHSS;
			Quinn-spaces)
			homeo\-morphism-in\-var\-iant geometric fundamental-classes exist.
			The construction combines three ideas:
			Firstly, instead of restricting geometric cycles by conditions on links only, a more flexible framework is built directly via
			geometric properties,
			secondly, controlled topology methods are used to give an accessible link-based criterion to detect suitable
			cycles
			and thirdly, a geometric argument is used to show, that these classes of cycles are
			suitable to study the transition to intrinsic stratifications.
			As an application, we give a construction of topologically (homeomorphism) invariant
			(homological) L-classes on MHSS Witt-spaces
			satisfying conditions on Whitehead-groups of links and the dimensional spacing
			of meeting strata.
			These L-classes agree, whenever those spaces are additionally pl-pseudo\-mani\-folds,
			with the Goresky--MacPherson L-classes.
		\end{abstract}
		
		
		
		\begin{keyword}
			Bordism\sep Generalized homology-theories\sep Singular spaces\sep Homotopically stratified spaces\sep
			Controlled topology\sep L-classes\sep Stratified transversality
			
			
			\MSC[2020]	
			55N20\sep	
			55N22\sep	
			55N33\sep	
			57N75\sep	
			57N80\sep	
			57R20		
			
		\end{keyword}
		
	\end{frontmatter}
	


	\newtheorem{Theorem}{Theorem}
	\newtheorem*{Theorem1}{Theorem \ref{thm:AlgTopCond}}
	\newtheorem*{Theorem2}{Theorem \ref{thm:WittBordism}}
	\newtheorem*{Lemma63}{Lemma \ref{lemma:SimCut}}

	\newtheorem{Lemma}{Lemma}[section]
	\newtheorem{Corollary}[Lemma]{Corollary}
	\newdefinition{Definition}[Lemma]{Definition}
	\newdefinition{Remark}[Lemma]{Remark}
	\newdefinition{Example}[Lemma]{Example}

	\newproof{Proof}{Proof}
	\newproof{ProofOfCorollary}{Proof of the corollary}
	\newproof{ProofOfThm1}{Proof of Theorem \ref{thm:AlgTopCond}}
	\newproof{ProofOfThm2}{Proof of Theorem \ref{thm:WittBordism}}

	\section{Introduction}
	\label{sec:introduction}
	
	As illustrated for example by Thom's \citep{ThomMexiko58} influential \citep{GoreskyMacPherson,SiegelWittSpaces}
	approach to the definition of L-classes
	through bordism-invariant signatures,
	geometric-cycle generalized homology-theories -- such as bordism(-type) \citep{QuinnBordismType,BRSMockBundles}
	or ad-theories \citep{LauresMcClureAdMulti,BanaglLauresMcClure} and their associated Quinn-spectra
	-- constitute a valuable tool in the study of characteristic classes and their transport-behavior.
	
	Geometric cycles carry information like signatures, when being represented by manifolds.
	Intersection-homology theory \citep{GoreskyMacPherson}
	provided a suitable bordism-invariant signature for singular-spaces as well:
	First \citep{GoreskyMacPherson} treated the case of spaces with only even-codimensional strata,
	then Siegel \citep{SiegelWittSpaces} introduced Witt-spaces generalizing these results
	and even further generalizations have been given by Banagl \citep{BanaglLClassOfNonWitt} using
	Lagrangian structures.
	However, the construction of a (generalized homology) bordism theory requires additionally suitable transversality
	results, to construct the inverse of the excision-isomorphism (or equivalently
	desuspension\Slash{}the Mayer--Vietoris boundary map).
	In the pl-category these transversality constructions can be carried out rather directly
	\citep{AkinBordism,SiegelWittSpaces,PardonIPBordism,FriedmanBordismPseudomanif,BanaglLauresMcClure},
	while in the topological world already the manifold-case \citep{KirbySiebenmannEssays,QuinnTransversality}
	is far from trivial.
	
	Summarized, this means, the topological manifold-case
	\citep{NovikovTopInv,KirbySiebenmannEssays,RanickiAlgebraicLTheory} and
	the singular pl-case \citep{ThomMexiko58,GoreskyMacPherson,SiegelWittSpaces,BanaglLauresMcClure}
	are quite well-understood. But neither are general topological manifolds triangulable, nor are
	pl-pseudomanifolds topological manifolds, which makes one wonder how to reconcile both worlds.
	We cannot fully resolve this question, but we can provide a tentative step towards understanding the unified picture,
	by giving inherently topological bordism-theories of singular spaces.
	We are working with Quinn's MHSS \citep{QuinnHSS} -- further requirements are indicated below -- which have been
	studied intensely for example by
	\citep{QuinnHSS, QuinnIntrinsicSkeleta,WeinbergerBook, CappellShanesonMappingCylinder, HTWW, HughesNbhEx02, FriedmanPD}. These spaces can be characterized as locally homotopy-conelike manifold-stratified spaces
	\citep{QuinnIntrinsicSkeleta}, in particular Siebenmann's CS-sets \citep{SiebenmannCS}
	and thus pseudomanifolds in the sense of \citep{GoreskyMacPherson} are MHSS.
	
	The treatment as presented here includes the topological manifold case as a special case, but the author
	wants to point out, that this is \emph{not} a claim of an alternative to Novikov's proof \citep{NovikovTopInv},
	but rather the manifold-transversality results \citep{KirbySiebenmannEssays,QuinnTransversality} which we implicitly rely on,
	build on the legacy of Novikov's ideas.
	Also our result is \emph{not} in strict generalization of the pl-results for example of \citep{BanaglLauresMcClure},
	because many pl-pseudomanifolds are not captured by our theory. This would require substantial improvements of
	the material presented in §\ref{sec:AlgTopCond}.
	Instead, the relevance of the approach presented in this paper is,
	that it allows the treatment of questions that cannot be
	answered independently for manifold and pl-case, like \inquotes{given a pl-pseudomanifold, are its L-classes
		invariant under \emph{topological} homeomorphisms?} or \inquotes{given a manifold and a (non-combinatorial) triangulation
		thereof which turns it into a pl-pseudomanifold, are Goresky--MacPherson L-classes the (Poincar\'{e}-duals of) the
		manifold L-classes?}.
	
	To this end, we give requirements on manifold homotopy stratified sets \citep{QuinnHSS}, that can be verified via a
	reasonably accessible algebraic topological condition on links,
	which in turn ensures the existence of geometric fundamental-classes
	in a generalized-homology bordism-theory.
	Namely, we require a condition on Whitehead groups related to the fundamental groups of links (see \inquotes{simple links},
	Def.\,\ref{def:Thm1Hyp}) that ensure the existence of controlled end-completions \citep{QuinnEndsOfMaps1}
	(a question closely related to h-cobordism triviality), which can be
	employed to obtain transversality results \citep{ConnollyVajiac}.
	Additionally a dimensional gap condition is required, which reflects the use of Whitney-tricks in end-completions and
	h-cobordism trivializations. This leads in §\ref{sec:AlgTopCond} to the result (where the transversality-construction
	required for excision is referred to as a \inquotes{cut}):

	\begin{Theorem1}
		Given a $5$-gapped MHSS $(X,\bdry X)$ with simple links,
		then $(X,\bdry X)$ has sufficiently many relative cuts.
	\end{Theorem1}
	
	The reason for requiring obstruction-\emph{groups} (rather than obstructions) to vanish
	and dimensional \emph{gaps} (rather than minimum dimensions) is in the symmetries required for the bordism-construction:
	Essentially, it must be possible to use excision\Slash{}desuspension arbitrarily (finitely) often
	(potentially while crossing with intervals in-between). So in contrast to \citep{ConnollyVajiac},
	we do not only make maps to $\Reals$ transverse to zero, but actually give a transversality-result for
	trivial normal-bundles of arbitrary rank.
	
	Transport under topological homeomorphisms, even unstratified ones, can be studied naturally in this setting:
	We use the intrinsic stratifications given by Quinn \citep{QuinnIntrinsicSkeleta} and their
	\emph{geometric} properties (in particular $|X \times \Reals| = |X|\times \Reals$ if $X^0 = \emptyset$) to show
	that the \emph{geometric} property of possessing \inquotes{sufficiently many cuts}
	required for excision is stable under taking cylinders in a sense suitable to give statements
	(see §\ref{sec:IntrinsicTransversality}) essentially of the form:
	\inquotes{If $X$ and its intrinsic stratification $|X|$ both have fundamental classes $[X], [|X|]$ in a (generalized homology)
		bordism theory, then there is another (generalized homology) bordism theory in which additionally a stratification
		of the cylinder $X \times I$ is allowed as a bordism $X \bordant |X|$.}
	From a practical perspective, requirements on $|X|$ (which we might not know)
	can be checked only in specific cases (see Example~\ref{example:SatisfyHypotheses}\,(d)).
	However, this means the question about topological transport can be answered via a
	condition on \emph{the underlying topological space}, independent of the choice of stratification,
	which is -- from a more philosophical perspective -- very appealing.
	In this spirit, the Theorem \ref{thm:AlgTopCond} indicated above \emph{detects} allowable cycles, thus spaces
	possessing fundamental classes, while the bordism-theory itself will be defined for cycles in the
	\inquotes{intrinsic completion} of §\ref{sec:IntrinsicTransversality}.
	The relation of stratifications and underlying spaces in the pl-setting was studied for example in
	\citep{FriedmanBordismPseudomanif}.
	
	Since the form of the additional bordisms (as stratifications of cylinders $X \times I$)
	is known, the entire construction will be seen to be compatible with the addition of further geometric rigidity,
	for example a Witt-condition \citep{SiegelWittSpaces,FriedmanPD} (as demonstrated in §\ref{sec:LClasses})
	leading to a signature-invariant, and thus L-classes -- the Goresky--MacPherson L-classes \citep{GoreskyMacPherson}
	in the case
	of pl-pseudo\-mani\-folds.
	
	Combining these three ingredients: (a) detection of allowable cycles and thus fundamental-classes
	via Theorem \ref{thm:AlgTopCond},
	(b) geometric \inquotes{completion} under cylinders to intrinsic stratifications (§\ref{sec:IntrinsicTransversality})
	and (c) compatiblility with Friedman's version of Witt-spaces \citep{FriedmanPD} for MHSS we obtain:
	
	\begin{Theorem2}
		(Shortened version, see §\ref{sec:DefsAndResults} and §\ref{sec:LClasses}):
		There is a generalized homology theory $\Bordism_*( \blank, \blank )$ on pairs of topological spaces,
		such that:
		\begin{enumerate}[(a)]
			\item Given a $5$-gapped, closed,
			oriented, MHSS Witt-space $X$ with simple links of dimension $n$ with
			dense top-stratum, there is a fundamental class
			$[X] := [X \xrightarrow{\id} X] \in \Bordism_n( X )$.
			\item These fundamental-classes are orientation-preserving
			stratified-homeo\-morphism  invariant.
			Under weak hypotheses, and if the underlying intrinsic stratification $|X|$ is again $5$-gapped with simple links, the homeomorphism need not be stratified.
			\item There is a group-homomorphism $\sigma_X : \Bordism_*(X) \rightarrow \Integers$,
			such that on these fundamental classes $\sigma_X([X]) = \sign(X)$ is the
			(middle perversity, middle dimension) rational intersection-homology pairing signature of the (MHSS) Witt-space $X$.
		\end{enumerate}
	\end{Theorem2}
	
	There are many relevant examples satisfying the hypotheses of this theorem,
	including topological manifolds, dimensionally gapped Witt-spaces with simply-connected holink-fibers or two-stratum,
	super-normal, dimensionally gapped Witt-spaces.
	For a summary of further examples see Examples \ref{example:SatisfyHypotheses}.
	The hypotheses are on \emph{pairs} of strata only, which is in line with the general
	philosophy of HSS.
	
	Many of the hypotheses required here may not be strictly necessary,
	but the approach presented provides an understanding of problems and necessary techniques --
	in particular from combining controlled topology \citep{QuinnEndsOfMaps2} with more geometric techniques
	\citep{HughesCylinders,ConnollyVajiac,QuinnIntrinsicSkeleta} -- involved in the study of
	existence and potentially \ParenthesizedHyphenation{non}uniqueness (see Rmk\,\ref{rmk:Choices}) of geometric
	representatives in generalized homology-theories.
	For a brief technical summary of the results presented in this paper, see also
	§\ref{sec:DefsAndResults} below.
	
	Finally, we want to point out, that there is a complementary, algebraic perspective
	\citep{RanickiAlgebraicLTheory,RanickiYamasakiControlledL},
	where similar (K-theoretic, see §\ref{sec:AlgTopCond}) difficulties are encountered.
	From this point of view, the difficulties that we find in constructing the bordism-theory
	instead appear when one tries to obtain fundamental classes, whose
	construction beyond the manifold-case \citep[§16]{RanickiAlgebraicLTheory}
	and pl-case \citep{BanaglLauresMcClure} again becomes difficult.
	See also \citep{RanickiWeiss10}.
	

	\section{Definitions and Statement of Results}\label{sec:DefsAndResults}
	
	The properties of a bordism theory, that is a generalized homology-theory where cycles are represented geometrically,
	depend vastly on what geometric representation is taken to mean precisely.
	Possible choices of geometric cycles are restricted by the symmetries
	of operations involved in constructing the generalized homology-theory.
	For the resulting theory to satisfy the excision-axiom,
	in particular some transversality properties are needed.
	
	In the pl-setting conditions on links (of the stratification) are a natural choice to restrict allowable cycles:
	Links behave relatively natural,
	because they are closely related to pl-links (see e.\,g.\ \citep{FriedmanBordismPseudomanif}).
	Such conditions on links are convenient, because they
	are automatically compatible with the symmetries of the cycles of a generalized bordism-theory.
	We will refer to link-based (and similar) conditions as auxiliary in §\ref{sec:CompatibleConditions},
	in the sense that they can be added without knowledge or modification of the internals of the
	transversality-construction used, thus are \inquotes{auxiliary} to the transversality statement.
	However this will not typically work in the \emph{topological} category: Already being a manifold is not auxiliary in
	this sense, as a \BiCollared{} closed subspace of a manifold need not itself be a manifold,
	a counterexample is Bing's \inquotes{dogbone}-space \citep{BingDogboneSpace}!
	For technical convenience we define boundaries of manifold-stratified spaces to be manifold-stratified, so this particular
	example\Slash{}problem is shifted to §\ref{sec:AlgTopCond}, where it is solved implicitly using controlled topology techniques.
	Nevertheless, other issues with links prevail (for example uniqueness of links, see \citep{HandelCS}),
	and while we do make heavy use of link-based conditions (see Theorem \ref{thm:AlgTopCond}),
	we will have to go beyond such \inquotes{auxiliary} conditions.
	
	To this end, we take a slightly different approach to \inquotes{compatible} conditions on cycles
	(compatible with the symmetries of generalized homology cycles, see above),
	via geometric properties satisfied by spaces.
	This, in particular turns out to be easier to reconcile with the intrinsic (stratification independent)
	topology of spaces. The resulting setup based on \inquotes{compatible conditions} is worked out in §\ref{sec:CompatibleConditions}
	and applied to bordism-theories in §\ref{sec:Bordism}, where Theorem \ref{thm:Bordism}, given such a compatible condition,
	determines the ensuing properties of an associated bordism-theory.
	We are working with singular (non-manifold spaces), and first fix some
	notation (see e.\,g.\ \citep{QuinnHSS,HTWW,HughesTeardrops}):
	
	\begin{Definition}\label{def:Basics}
		A \defName{stratified space} $X$ is a topological space together with a filtration by closed	subsets
		\begin{equation*}
			X = X^n \supset X^{n-1} \supset \ldots \supset X^0 \supset X^{-1} = \emptyset
		\end{equation*}
		with $X^k$ called the $k$-skeleton of $X$, and $X_k := X^k-X^{k-1}$ the $k$-stratum.
		We allow strata to be disconnected (for indexing them by dimension in the manifold-strata case),
		but assume there is only a finite number of components.
		
		A stratified pair of stratified spaces $(X,\bdry X)$ is \defName{(weakly) manifold-stratified}, if
		\begin{enumerate}[(i)]
			\item $\forall k$ the $k$-stratum $X_k$ is a (topological) $k$-manifold (potentially with boundary).
			\item the \defName{boundary} $\bdry X$ of $X$ is the union $\cup_k \bdry X_k$
			of the manifold-boundaries of the strata $X_k$.
			\item the boundary is \defName{stratified collared}, i.\,e.\ $\bdry X \subset X$ is closed and there is
			a stratified (mapping strata to strata; the left-hand-side is given the product stratification)
			map $c : \bdry X \times [0,\infty) \rightarrow X$,
			which is a homeomorphism to its image $\im(c)$ and $\im(c)$ is an open neighborhood of $\bdry X$
			in $X$.
		\end{enumerate}
		and we drop the "weakly" if additionally (cf.\ Rmk.\,\ref{rmk:FrontierCondtionForMHSS})
		\begin{enumerate}[(i)]
			\setcounter{enumi}{3}
			\item $X$ and $\bdry X$ satisfy the local frontier-condition: For any two
				path-connect\-ed components $A$, $B$ of strata $X_i, X_k$:
				$A \cap \closure(B) \neq \emptyset \Rightarrow A \subset \closure(B)$
		\end{enumerate}
		
		A separated, metric weakly manifold-stratified space $(X,\bdry X)$ is a \defName{manifold homotopy stratified set (MHSS)},
		if for every pair of \defName{meeting (components of) strata}, that is $(S,T)$ with $\closure(S) \cap T \neq \emptyset$ 
			and $S\neq T$, satisfies:
		\begin{enumerate}[(i)]
			\item $T \subset S \cup T$ is \defName{tame}, i.\,e.\ it is a nearly stratum-preserving
			strict neighborhood deformation retract. This means,
			there is a retraction $R_t$ such that $R_t|_T = \id_T$
			and $\forall t<1 : R_t(S)\subset S$, while $\im(R_1) \subset T$.
			\item the \defName{homotopy-link} $\holink(S \cup T, T) := \{ \gamma : [0,1] \rightarrow S \cup T |
			\gamma(0) \in T$ and $\forall t > 0 : \gamma(t) \in S \}$ is such that evaluation at zero,
			$p: \holink(S \cup T, T) \rightarrow T, \gamma \mapsto \gamma(0)$ is a fibration.
		\end{enumerate}
		Some authors include the hypothesis of \inquotes{compactly dominated local holinks} \citep{HughesTeardrops}
		in the definition of MHSS.
		Here, this choice only matters for Witt-spaces in §\ref{sec:LClasses}.
		This is a local condition and satisfied on the compact spaces used in §\ref{sec:LClasses},
		by \citep[Lemma\,5.3 (p.\,318)]{HughesTeardrops}, thus we may safely add it in hindsight there.
	\end{Definition}
	
	\begin{Remark}\label{rmk:FrontierCondtionForMHSS}
		MHSS as defined above satisfy the local frontier-condition automatically, thus are always
		(non-weakly) manifold-stratified (see e.\,g.\ \citep[discussion before Prop.\,10.2 (p.\,30)]{HughesStratPath}):
		If $A \cap \closure(B) \neq \emptyset$, then tameness implies the $\holink(A\cup B,A) \neq \emptyset$.
		If there were $x\in A - \closure(B)$, using path-connectedness of $A$, we obtain a path $\gamma$ from a point
		$y \in A \cap \closure(B) \neq \emptyset$ to $x$, produce a lift $\bar\gamma$ in the holink, thus
		$\bar\gamma(t)((0,1]) \subset B$, hence $\bar\gamma(t)(1) \in \closure(B)$, in particular
		$x = \bar\gamma(1)(1) \in \closure(B)$.
	\end{Remark}
	
	The difficulty in building a geometric-cycle generalized homology-theory is in the transversality statement
	needed for the excision
	axiom. This can be formalized as (the compatible version, see §\ref{sec:CompatibleConditions}, of) a \inquotes{cutting} property:
	
	\begin{Definition}\label{def:Cut}
		Given a manifold-stratified space $(X,\bdry X)$,
		a \defName{cut} of $X$ is a triple $(X_\geq, X_\leq, X_{\{0\}})$
		of manifold-stratified sub-spaces
		(subsets with the induced stratification by intersection with the filtration),
		with $X_\geq \cup X_\leq = X$ and $X_\geq \cap X_\leq = X_{\{0\}}$, where the
		(collared) boundaries are $\bdry X_\geq = (\bdry X)_\geq \cup X_{\{0\}}$ with
		$(\bdry X)_\geq := X_\geq \cap \bdry X$, similarly for $X_\leq$,
		and $\bdry X_{\{0\}} = (\bdry X)_{\{0\}} := X_{\{0\}} \cap \bdry X$.
		
		Given closed disjoint subsets $B_\pm \subset X$,
		a cut ($X_\geq$, $X_\leq$, $X_{\{0\}}$) is \defName{between $B_+$ and $B_-$},
		if $X_\geq \cap B_- = \emptyset$ and $X_\leq \cap B_+ = \emptyset$.
		
		We call a cut $(X_\geq, X_\leq, X_{\{0\}})$ of $X$ \defName{boundary-compatible},
		if there is a boundary collar $b$ of $\bdry X \subset X$ and a cut
		$((\bdry X)_\geq, (\bdry X)_\leq, (\bdry X)_{\{0\}})$ of $\bdry X$,
		such that $b( ((\bdry X)_\geq, (\bdry X)_\leq, (\bdry X)_{\{0\}}) \times [0, \infty))$
		agrees with the cut of $X$ on an open neighborhood of $\bdry X$.
		
		A manifold-stratified space $X$ has \defName{sufficiently many absolute cuts} if given
		closed disjoint subsets $B_\pm \subset X$
		there is a boundary-compatible cut $(X_\geq$, $X_\leq$, $X_{\{0\}})$ of $X$ between $B_\pm$.
		
		A manifold-stratified space $X$ has \defName{sufficiently many relative cuts} if given
		closed disjoint subsets $B_\pm \subset X$	
		and a boundary-compatible cut $(D_\geq$, $D_\leq$, $D_{\{0\}})$ of $D := X^k \cup C$ between $D \cap B_{\pm}$,
		where $C \subset X$ is an open subset of $X$,
		then given any open neighborhood $U$ of the complement $X - C$,
		there is a boundary-compatible cut $(X_\geq, X_\leq, X_{\{0\}})$ of $X$ between $B_\pm$,
		such that on $D' := X^k \cup (X-U)$ the restrictions
		$(D' \cap X_\geq, D' \cap X_\leq, D' \cap X_{\{0\}}) = (D' \cap D_\geq, D' \cap D_\leq, D' \cap D_{\{0\}})$
		agree,
		and the boundary-collars of $X_\geq$ and $X_\leq$ extend
		those of $X^k$ as germs (i.\,e.\ the new collar is allowed to be thinner; 
		for an explanation, why the collars are extended only from $X^k$ not from all of $D$,
		see Rmk.\,\ref{rmk:CutUniqueness} and the proof of Theorem \ref{thm:AlgTopCond} in §\ref{sec:AlgTopCond}).	
	\end{Definition}
	
	Using this notion, being compatible as opposed to auxiliary (see above)
	will be understood as producing the correct cuts rather than being
	automatically consistent with arbitrary cuts.
	
	For example manifolds have sufficiently many relative cuts \citep{KirbySiebenmannEssays,QuinnTransversality}:
	If $X=M$ is a manifold, let $g: M \rightarrow [-1,1]$ a continuous map with $g|_{B_{\pm}} = \pm 1$,
	then make $g$ transverse to $0$ rel $g^{-1}(\{\pm 1\})$ and let $X_\geq = g^{-1}([0,1])$, $X_\leq = g^{-1}([-1,0])$
	and $X_{\{0\}} = g^{-1}(\{0\})$. This is a cut between $B_{\pm}$, and can be obtained rel open subset away from the
	boundary and\Slash{}or of the boundary (cut the boundary, extend as a product into the boundary collar,
	and apply manifold-transversality rel the outer half of the collar).
	
	A more interesting class of examples is described by the following
	sufficient conditions obtained through controlled-topology methods \citep{QuinnEndsOfMaps2} and results of \citep{ConnollyVajiac} in §\ref{sec:AlgTopCond}
	for MHSS \citep{QuinnHSS}.
	
	\begin{Definition}\label{def:Thm1Hyp}
		A MHSS $X$ is \defName{$d$-gapped} if meeting strata (Def.\,\ref{def:Basics}) differ in dimension by at least $d$.
		
		A MHSS $X$ has \defName{simple links}, if
		each component of a fiber of a holink
		has fundamental group $\pi$ such that $\Wh(\pi \times \Integers^k) = 0$ for all $k\geq 0$.
	\end{Definition}
	
	\begin{Theorem}\label{thm:AlgTopCond}
		Given a $5$-gapped MHSS $(X,\bdry X)$ with simple links,
		then $(X,\bdry X)$ has sufficiently many relative cuts.
	\end{Theorem}
	
	Note, that the hypothesis of this theorem depends only on codimensions and the homotopy-type of holink-fibers, thus can
	rather easily be seen to be auxiliary in the sense of §\ref{sec:CompatibleConditions}, so these spaces
	actually can be cut repeatedly and are thus compatible in the sense of §\ref{sec:CompatibleConditions}.
	
	In §\ref{sec:IntrinsicTransversality} transversality is seen
	to be \inquotes{almost} intrinsic, more precisely: To occur
	simultaneously on intrinsic stratifications for spaces with conelike $4$-skeleton either
	away from certain \inquotes{strange points} (e.\,g.\ if $X^0 = \emptyset$) or stably (see the Lemma below).
	
	\begin{Definition}\label{def:StrangePoints}
		Let $|X|$ denote the \defName{intrinsic stratification} of a MHSS $X$ of \citep{QuinnIntrinsicSkeleta},
		see §\ref{sec:IntrinsicTransversality} for details.
		
		A MHSS $X$ has \defName{no strange points},
		if $|X \times \Reals| = |X| \times \Reals$.
		Equivalently 
		(see Rmk.\,\ref{rmk:StrangePoints} and \citep{QuinnIntrinsicSkeleta}),
		$X^0$ must not contain points whose homotopy-link-fibers in any meeting strata $|X|_k$
		of $|X|$ have the homology of a ($k-1$)-sphere while being non-simply-connected.
	\end{Definition}
	
	So \inquotes{strange points} are really represented by points, and do not exist if $X^0 = \emptyset$
	(see Rmk.\,\ref{rmk:StrangePoints}).
	In low dimensions geometric topology methods cease to work, which is where the $5$-gapped hypothesis comes from,
	but which will also lead to some local
	conelikeness requirements (on stratifications, as this is not well-defined on underlying spaces
	\citep{HandelCS}). This concerns (at most) $4$-skeleta, which for $5$-gapped spaces as used
	in Theorem \ref{thm:AlgTopCond} are disjoint unions of manifolds, thus trivially locally conelike.
	
	\begin{Definition}
		A stratified space $X$ is \defName{locally conelike} \citep{SiebenmannCS,QuinnIntrinsicSkeleta},
		if for each $x \in X_j$, there is a compact stratified space $L$ and an open neighborhood $U$ of $x$
		together with a stratified homeomorphism	$\varphi : \opencone{L} \times \Reals^j \rightarrow U$
		with $\varphi( v, 0 ) = x$. The $k$-skeleton of the \defName{stratified cone}
		is the image of $L^{k-1} \times [0,1)$ under the quotient-map $L \times [0,1)  \rightarrow \opencone(L)$
		if $k > 0$ and the cone-vertex $\{v\}$ if $k=0$. The left-hand-side is then given the product-stratification,
		the right-hand-side the induced (subspace) one from $X$.
	\end{Definition}
	
	Cuts are then simultaneous in a sense illustrated for example by:
	
	\begin{Lemma63}
		Given a MHSS $X$ with locally conelike $3$-skeleton and
		a bound\-ary-compatible cut $(X_\geq,X_\leq,X_{\{0\}})$ of $X$,
		the triple $(|X_\geq \times \Reals|,|X_\leq \times \Reals|,|X_{\{0\}} \times \Reals|)$
		is a boundary-compatible cut of $|X \times \Reals|$.
	\end{Lemma63}
	
	This means for the resulting bordism-theories,
	that sufficient conditions to satisfy the excision axiom can in some cases be given independently of the
	choice of stratification, for example rendering fundamental classes topologically invariant for
	suitable spaces. Somewhat surprisingly, for the examples provided by Theorem \ref{thm:AlgTopCond},
	this invariance can be detected from the intrinsic stratification, thus the underlying space
	allows to conclude when fundamental-classes are stratification-independent.
	
	Finally, the resulting theory will be applied with a Witt-condition \citep{SiegelWittSpaces,FriedmanPD},
	to show that on suitable (see Example \ref{example:SatisfyHypotheses}) MHSS
	homological L-classes can be defined, that are, for suitable spaces (see below)
	invariant under orientation-preserving \emph{unstratified} and \emph{topological} homeomorphisms.
	
	\begin{Theorem}\label{thm:WittBordism}
		There is a generalized homology theory defined on pairs of topological spaces,
		realized as oriented bordism of certain\footnote{The geometric cycles used are
			certain \inquotes{intrinsic multi-gluings}, see §\ref{sec:IntrinsicTransversality}, of $5$-gapped
			MHSS with simple links, which are additionally Witt, with dense top-stratum.
			The proof of this theorem in §\ref{sec:LClasses} provides the details on this.}
		MHSS Witt-spaces with dense top-stratum
		$\Bordism_*( \blank, \blank )$,	such that:
		\begin{enumerate}[(1)]
			\item The graded group $\Bordism_*( Z, A )$ is a module over the oriented mani\-fold-bordism ring
			$\Omega_*^{\txt{STOP}}( \ptSpace )$.
			\item Given a $5$-gapped, closed,
			oriented, MHSS Witt-space $X$ of dimension $n$ with simple links with
			dense top-stratum, there is a fundamental class
			$[X] := [X \xrightarrow{\id} X] \in \Bordism_n( X )$.
			\item If $h: X\rightarrow Y$ is a stratified orientation-preserving homeomorphism,
			then $h_*[X] = [Y]$ for these fundamental classes (if they exist).
			If additionally both $X$ and $Y$ are stratified without strange points,
			and the underlying intrinsic stratification $|X|$ of the source $X$
			is again $5$-gapped with simple links, the hypothesis on $h$ of being stratified
			can be dropped.
			\item There is a group-homomorphism $\sigma : \Bordism_*(Z) \rightarrow \Integers$,
			such that on these fundamental classes $\sigma([X]) = \sign(X)$ is the
			(middle perversity, middle dimension) rational intersection-homology pairing signature of the (MHSS) Witt-space $X$.
		\end{enumerate}
	\end{Theorem}
	
	Actually, any space allowed as cycle has a fundamental class, so this immediately generalizes
	to spaces with boundary and fundamental classes in $\Bordism_n( X, \bdry X )$,
	but hypotheses start to get rather hard to read.
	
	\begin{Example}\label{example:SatisfyHypotheses}
		The hypotheses can be checked as follows:
		\begin{enumerate}[(a)]
			\item Any topological closed oriented $n$-manifold
			is an oriented $5$-gapped MHSS Witt-space of dimension $n$ with simple links, no strange points and
			dense top-stratum, by manifold-transversality and Poincar\'{e}-dual\-ity.
			Technically this uses, that manifold-boundaries are collared \citep{BrownCollars}.
			Further $|M| = M$, so the fundamental class of a manifold is (unstratified) homeomorphism invariant
			(as long as any given other stratification has a fundamental class and no strange points).
			
			\item A MHSS $X$ -- for example a CS-set
			(see e.\,g.\ \citep[Thm.\,2 (p.\,239) and \inquotes{Converse} (p.\,240)]{QuinnHSS}),
			including pseudomanifolds --
			with holink-fibers having free abelian fundamental-groups (for example simply-connected links,
			including 2-stratum \inquotes{super-normal} spaces, see e.\,g.\ 
			\citep[§12.1 (p.\,202f)]{WeinbergerBook};
			the restriction to 2-stratum spaces accounts for the super-normal condition usually being given
			on stratified local holinks, whose strata are the links used here)
			has simple links by Bass--Heller--Swan.
			
			\item A MHSS $X$ with $X^0 = \emptyset = (\bdry X)^0$
			has no strange point (see Rmk.\,\ref{rmk:StrangePoints}).			
			For the calculation of L-classes, spaces will be \inquotes{stabilized} by products with high-dimensional spheres, so strange points never occur.
			
			\item $|X|$ being $5$-gapped and having simple links can sometimes be checked from $X$.
			For example if a $5$-gapped MHSS $X$ with simple links has at most $2$ meeting strata,
			then $|X|$ also is $5$-gapped with simple links, see Example \ref{example:Intrinsic5GappedSimpleLinks} (a).
			A similar argument applies if $|X|$ is a topological manifold,
			as in the case, when we want to compare Goresky--MacPherson L-classes
			of a pl-pseudomanifold stratification of a topological manifold,
			with the manifold L-classes,
			see Example \ref{example:Intrinsic5GappedSimpleLinks} (b).
			
			\item A $5$-gapped MHSS $X$ is Witt if and only if it satisfies the usual link condition (using the upper-middle perversity,
			the middle dimension inter\-sect\-ion-homology of the stratified local holink vanishes), by \citep{FriedmanPD}.
			Thus if $X$ is a pseudomanifold, it is 
			Witt as a MHSS if and only if it is Witt as pseudomanifold.
		\end{enumerate}
	\end{Example}
	
	\begin{Corollary}
		If $X$ is a $5$-gapped closed, oriented MHSS Witt-space
		with simple links and dense top-stratum,
		then there are (homological) L-classes $L_*(X)$, which agree with the
		Goresky--MacPherson L-classes if $X$ is a pl-pseudomanifold.
		
		These L-classes are stratified homeomorphism invariant
		and unstratified homeomorphism invariant if $X, Y$ are additionally without strange points
		and $|X|$ is $5$-gapped with simple links,
		that is, given a (stratified in the first case) orientation-preserving homeomorphism $h: X \rightarrow Y$, then 
		\begin{equation*}
			h_*(L_*(X)) = L_*(Y)\txt{.}
		\end{equation*}
		In particular, Goresky--MacPherson L-classes of
		$5$-gapped, closed, oriented pl-pseudomanifold Witt-spaces with simple links
		are stratified orientation-pre\-serv\-ing homeomorphism-invariant and
		unstratified orientation-preserving homeomorphism-invariant, if $X, Y$ are additionally without strange points
		and $|X|$ is $5$-gapped with simple links.
	\end{Corollary}
	
	Invariance of Goresky--MacPherson L-classes
	under pl-homeomorphisms has been treated in detail by \citep{BanaglLauresMcClure}.
	The (topological) homeomorphism case was given for smooth manifolds by Novikov \citep{NovikovTopInv},
	and in the stratified case has been discussed before \citep{CSW91} and \citep[p.\,209f]{WeinbergerBook},
	but the details and geometry behind the problem seem much less understood than in the pl-case.
	
	There are $5$-gapped MHSS with simple links,
	that are not pl-pseudo\-mani\-folds (for example non-triangulable topological manifolds),
	but of course the converse is also true, so the result here is neither strictly weaker nor strictly stronger than
	the known pl-results.
	
	The main technical results of this paper are given in §\ref{sec:AlgTopCond},
	giving the proof of Theorem \ref{thm:AlgTopCond},
	and §\ref{sec:IntrinsicTransversality}, giving constructions to reconcile the theory with intrinsic stratifications.
	The sections §\ref{sec:CompatibleConditions} and §\ref{sec:Bordism} provide the technical definitions for
	the remainder of the paper and the construction of bordism-theories in this context.
	The final §\ref{sec:LClasses} applies this to singular L-classes, in particular proving Theorem \ref{thm:WittBordism}.

	\section{Compatible Conditions}\label{sec:CompatibleConditions}
	
	To combine transversality with the symmetry-requirements -- persistence to cuts, under products and with gluings --
	needed for bordism-theories, we formulate conditions on properties\Slash{}classes of spaces.
	One could probably formulate this on suitable categories, however neither \inquotes{bordism-categories} nor
	\inquotes{bordism-type theories}\Slash{}
	\inquotes{ad-theories} are really what we want (see Rmk.\,\ref{rmk:CompatibleCategories}).
	To avoid huge technical overhead, we use the following more direct definitions.
	This means that orientation and compactness have to be added
	\inquotes{ad hoc} in §\ref{sec:Bordism}, but does not severely affect the theory.
	
	\begin{Definition}\label{def:CProperties}
		Let $\mathcal{C}$ a predicate of MHSS (i.\,e.\ $\mathcal{C}(X)$ is either true or false
		for a given MHSS $X$).
		
		\begin{enumerate}[(i)]	
			\item	A stratified space $(X,\bdry X)$ has \defName{$\mathcal{C}$-boundary}, if its boundary satisfies
			$\mathcal{C}(\bdry X)$.
			
			\item	A cut $(X_\geq, X_\leq, X_{\{0\}})$ of $X$ is a \defName{$\mathcal{C}$-cut}, if
			$X_\geq$ and $X_\leq$ have $\mathcal{C}$-boundary.
			
			\item	A stratified space $X$ has \defName{sufficiently many (absolute) $\mathcal{C}$-cuts} if
			given closed disjoint subsets $B_+, B_- \subset X$
			there is a $\mathcal{C}$-cut beetween $B_\pm$.
			
			\item	A stratified space $X$ with $\mathcal{C}(X)$ has the \defName{$\mathcal{C}$-product property} if
			given an arbitrary (oriented) manifold $M$ also $\mathcal{C}(X \times M)$ holds.
		\end{enumerate}
		
		Next, let $\mathcal{C}$ be given, such that for MHSS $X$, $Y$ with $\dim(X) = \dim(Y)$
		and such that given $Z \subset \bdry X, \bdry Y$ with
		$Z, \closure^{\bdry X}(\bdry X - Z), \closure^{\bdry Y}(\bdry Y - Z)$ all of $\dim(X)-1 = \dim(Y)-1$
		or empty,
		such that $\mathcal{C}$ is satisfied, then $\mathcal{C}(X \cup_Z Y)$ is satisfied
		(here, $Z = \emptyset$ is allowed, so this includes disjoint unions).
		And further, $\mathcal{C}(X)$ implies that $X$ has the $\mathcal{C}$-product property,
		then for all $X$ with $\mathcal{C}(X)$ satisfied
		\begin{enumerate}[(a)]
			\item if $X$ having $\mathcal{C}$-boundary implies $X$ has sufficiently many $\mathcal{C}$-cuts,
			then we call $\mathcal{C}$ \defName{compatible}.
			\item if $X$ has (automatically) $\mathcal{C}$-boundary, then we call $\mathcal{C}$ \defName{auxiliary}.
		\end{enumerate}
	\end{Definition}
	
	Some examples are given in Lemma \ref{lemma:ExamplesOfCompatibleConditions}.
	A $\mathcal{D}$ satisfying (b) is \inquotes{auxiliary} in the sense that choice of cut of $X$ with $\mathcal{D}(X)$ is
	automatically a $\mathcal{D}$-cut.
	An important example of compatible conditions
	will be provided by Theorem \ref{thm:AlgTopCond}
	together with the observation, that its hypothesis is auxiliary (cf.\ the next
	Lemma \ref{lemma:IntercompatibilityOfConditions}).
	In particular, non-trivial examples exist and can be detected in a reasonable way.
	
	\begin{Lemma}\label{lemma:IntercompatibilityOfConditions}
		If $\mathcal{C}$ is compatible and $\mathcal{D, D'}$ are auxiliary,
		then
		\begin{enumerate}[(i)]
			\item $(\mathcal{D} \txt{ and } \mathcal{D'})$ is auxiliary, and so is $(\mathcal{D} \txt{ or } \mathcal{D'})$.
			\item $(\mathcal{C} \txt{ and } \mathcal{D})$ is compatible.
			\item if $\mathcal{D}(X)$ $\Rightarrow$ $X$ has sufficiently many absolute cuts,
			then $\mathcal{D}$ is compatible.
		\end{enumerate}
	\end{Lemma}
	
	\begin{Proof}
		Both (i) and (ii) are immediate from the definitions.	
		Part (iii) follows directly from cuts existing and again satisfying $\mathcal{D}$.
	\end{Proof}
	
	Our main examples for such conditions are:
	
	\begin{Lemma}\label{lemma:ExamplesOfCompatibleConditions}
		The following conditions are auxiliary in this sense:
		\begin{enumerate}[(1)]
			\item $X$ has dense top-stratum, $\closure(X-X^{n-1}) = X$
			\item any requirement which depends only on
			codimension and homotopy-type of
			holink-fibers.
			This is also true for the stratified homotopy-type of stratified holinks, e.\,g.\ for
			the \inquotes{local homotopy-links} in the Witt-condition \citep{FriedmanPD}, see §\ref{sec:LClasses}.
			We will not use stratified holinks elsewhere.
			\item $X$ has no codimension $1$ stratum, $X^{n-1} = X^{n-2}$
			\item $X$ satisfies the hypothesis of Theorem \ref{thm:AlgTopCond}
		\end{enumerate}
		The following conditions are compatible in this sense:
		\begin{enumerate}[(i)]
			\item $X$ is $5$-gapped with simple links (by Theorem \ref{thm:AlgTopCond} proved in §\ref{sec:AlgTopCond}).
			\item The condition given in Lemma \ref{lemma:BordismSpaceIntrinsic}.
		\end{enumerate}
	\end{Lemma}
	
	\begin{Proof}
		First note, that any local condition (a condition that is satisfied on $X$ if it is on open covers)
		automatically is preserved to gluings as in Def.\,\ref{def:CProperties}.
		This applies to all conditions given above except (ii).
		
		The auxiliary conditions:
		(1) is immediate from boundaries being collared,
		(3) and (4) follow from (2), so we only need to show (2).
		
		Locality of (2) is clear (see the $\epsilon$-holinks in \citep{QuinnHSS}),
		it descends to collared boundaries, because the stratified
		collar induces a fibered deformation-retraction of the (stratified and pairwise) holinks.
		Further
		$X\times M$ has the same holinks and codimensions as $X$ and is a MHSS \citep{QuinnHSS}.
		
		The compatible conditions:
		(i) will follow from Theorem \ref{thm:AlgTopCond} together with part (4) above and
		Lemma \ref{lemma:IntercompatibilityOfConditions}\,iii, see also Cor.\,\ref{cor:AlgTopCond}.
		(ii) will be proved in §\ref{sec:IntrinsicTransversality}.
	\end{Proof}

	\section{Bordism Theories}\label{sec:Bordism}
	
	Most of the construction can be done in a standard way.
	The main technical difficulty is in the excision-axiom (equivalently,
	the (de)suspension-iso\-morph\-ism\Slash{}Mayer--Vietoris for the
	reduced theory), which requires a trans\-versality argument.
	
	\begin{Lemma}\label{lemma:BordismDef}
		Given a compatible or auxiliary condition $\mathcal{C}$,
		there are \defName{bordism groups in degree $n$}
		with elements
		\begin{align*}
			&\Omega^{\mathcal{C}}_n( Y, Z ) :=
			\{ f : (X,\bdry X) \rightarrow (Y,Z) \txt{ continuous } | \txt{ either } X=\emptyset	\\
			&\quad\txt{or }\dim(X) = n \txt{, }
			\mathcal{C}(X)\txt{, }\mathcal{C}(\bdry X)  \txt{ and } X \txt{ is a compact MHSS }
			\} / \bordant
		\end{align*}
		where the relation $(f : (X,\bdry X) \rightarrow (Y,Z)) \bordant (g : (X',\bdry X') \rightarrow (Y,Z))$ means, there is
		a compact $(W,\bdry W)$ with $\mathcal{C}(W)$ and $\bdry W = \bdry_0 W \cup \bdry_Z W \cup \bdry_1 W$
		with $\mathcal{C}(\bdry_0 W), \mathcal{C}(\bdry_1 W), \mathcal{C}(\bdry_Z W)$ such that
		$\bdry_0 W \cap \bdry_1 W = \emptyset$
		and $\bdry_0 W \cap \bdry_Z W = \bdry (\bdry_0 W)$,
		similarly for $\bdry_1 W$
		and $\bdry (\bdry_Z W) = \bdry (\bdry_0 W) \cup \bdry (\bdry_1 W)$
		where also $\mathcal{C}(\bdry (\bdry_0 W))$ and $\mathcal{C}(\bdry (\bdry_1 W))$
		and $\bdry(\bdry (\bdry_0 W)) = \emptyset = \bdry(\bdry (\bdry_0 W))$,
		together with $F : ( W, \bdry_Z W ) \rightarrow ( Y, Z )$ such that $\bdry_0 W = X$, $F|_{\bdry_0 W} = f$ and
		$\bdry_1 W = -Y$ (with reversed orientation in the oriented case, $\bdry W$ thus $\bdry_* W$ is given the
		induced orientation from $W$), $F|_{\bdry_1 W} = g$.
		
		The group-operation is given by disjoint union. Inverses are obtained by reversing orientation
		(or by doing nothing in the unoriented case, where everything is two-torsion).
		The neutral element (in any degree) is the empty set $\emptyset_n$ (formally considered as a space of dimension $n$).
		
		This is a functor via the induced maps given by composition: $\varphi_*([f]) := [\varphi \circ f]$.
		Finally, the graded abelian group $\Omega^{\mathcal{C}}_*( Y, Z )$ obtained by formally summing all degrees,
		is a module over the coefficient-ring of (oriented or unoriented respectively) topological manifold bordism.
	\end{Lemma}
	
	\begin{Proof}
		Bordism is an equivalence relation:
		It is reflective, because the product with the manifold $I$ again satisfies $\mathcal{C}$,
		$\mathcal{C}(X) \Rightarrow \mathcal{C}(X \times I)$ and $F = f \pi_X$ is continuous
		(this also shows the claim about inverses),
		it is symmetric by reversing the orientation on the bordism $W$ (or trivially so in the unoriented case)
		and it is transitive by gluing bordisms, which is allowed (see Def.\,\ref{def:CProperties})
		because boundary-components were required to satisfy $\mathcal{C}$.
		
		The functoriality is clear (composition of induced maps is composition of maps on representatives),
		and since $\mathcal{C}(X) \Rightarrow \mathcal{C}(X \times M)$ (for oriented $M$ if $\mathcal{C}$ implies
		orientability),
		scalar-multiples by manifolds over the point $f_M : M \rightarrow \{ \txt{pt} \}$ are given by
		cartesian products $[f_M] [f] := [f \times f_M]$.
	\end{Proof}
	
	\begin{Lemma}\label{lemma:BordismProperties}
		The boundaries $\bdry_n : \Omega^{\mathcal{C}}_n( Y, Z ) \rightarrow \Omega^{\mathcal{C}}_{n-1}( Z )$ given by
		$\bdry([f : (X,\bdry X) \rightarrow (Y, Z)]) := [f|_{\bdry X}]$ are natural transformations,
		such that $(\Omega^{\mathcal{C}}_*,\bdry_*)$ satisfies all axioms of a generalized homology
		\emph{except} potentially excision.
	\end{Lemma}
	
	\begin{Proof}
		Note, that if $F$ on $W$ is a bordism $f \bordant g$, then
		$F|_{\bdry_Z W}$ is a bordism $\bdry(f) \bordant \bdry(g)$, so the construction is well-defined.
		Since the group-structure is given by disjoint union and clearly $\bdry(f \sqcup g ) = \bdry(f) \sqcup \bdry(g)$,
		these are group-homomorphisms.
		They are natural transformations, because restriction of maps commutes with composition of maps.
		
		Homotopy-Axiom: If $f \homot_H g$, then $H : X \times I \rightarrow Y$ is a bordism by
		$\mathcal{C}(X) \Rightarrow \mathcal{C}(X \times I)$.
		Additivity-Axiom: Clearly the bordism-groups of disjoint unions are subgroups of the product, by compactness
		elements have finitely many entries.
		Exactness-Axiom: This part is a bit more lengthy, but elementary, so details are omitted. Only gluing (see transitivity),
		and some standard-arguments are needed.
		One needs to use bordisms as cycles and cycles as bordisms,
		so exactness forces us to use the same condition $\mathcal{C}$ both on geometric cycles and bordisms.
	\end{Proof}
	
	\begin{Lemma}\label{lemma:Excision}
		If $\mathcal{C}$ is compatible, then $(\Omega^{\mathcal{C}}_*,\bdry_*)$ also
		satisfies excision, thus is a generalized homology-theory.
	\end{Lemma}
	
	\begin{Proof}
		Given $(Y,Z)$ and $\closure(B) \subset \interior(Z)$, we need to show
		$\incl_* : \Omega^{\mathcal{C}}_k( Y-B, Z-B ) \rightarrow \Omega^{\mathcal{C}}_k( Y, Z )$ is an isomorphism.
		To this end, we construct an inverse $\psi$.
		Given $f : (X,\bdry X) \rightarrow (Y, Z)$, let $B_+ := X - \interior( f^{-1}(Z) )$
		and $B_- := f^{-1}(\closure(B))$. These are closed and disjoint by $\closure(B) \subset \interior(Z)$.
		By $\mathcal{C}$ being compatible, there is a $\mathcal{C}$-cut $(X_\geq, X_\leq, X_{\{0\}})$
		between $B_+$ and $B_-$, define $\psi([f]) := [f|_{X_\geq}]$.
		
		This is well-defined:
		(i) By choice of $B_-$ and $X_\geq \cap B_- = \emptyset$ (by definition of the cut) $f|_{X_\geq}$ maps to $Y-B$.
		(ii) Since this is a $\mathcal{C}$-cut , $\mathcal{C}(X_\geq)$ is satisfied.
		(iii) It is independent of the choice of the cut: Given a different cut $X'_\geq$, we may assume $X'_\geq \subset X_\geq$,
		otherwise find a third cut $X''_\geq$ between $X'_\geq \cup X_\geq$ and $B_-$ and use transitivity.
		Then $X_\geq \bordant X'_\geq$ via $W := X_\geq \times I$ and $f|_{X_\geq} \times \id_I$, using
		$\bdry_0 W = X_\geq \times \{0\}$, $\bdry_1 W = X'_\geq \times \{1\}$ and
		$\bdry_Z W = \closure(X_\geq - X'_\geq) \times \{1\} \cup \bdry X_\geq \times I$,
		which is allowed\Slash{}gets mapped to $Z$ by the choice of $B_+$ and all cuts being $\mathcal{C}$-cuts.
		(iv) It is independent of the choice of different representatives, by cutting the bordism $W$ (and restricting $F$)
		rel the cuts in the boundary-components (the representatives) which yields a bordism,
		in $\Omega^{\mathcal{C}}_k( Y-B, Z-B )$, through $W_\geq$.	
		
		Further $\psi$ is a group-homomorphism (the cut of a disjoint union, is a disjoint union of cuts),
		and $\psi \circ \incl_* = \id$, because $\psi$ is independent of the choice of representative,
		and if given $\incl_*([f])$, a possible choice is $X_\geq = X$.
		Finally, it is left to check that $\incl_* \circ \psi = \id$, which can be done by the same trick used above:
		Define $W := X \times I$, with $\bdry_0 W = X \times \{0\}$, $\bdry_1 W = X_\geq \times \{1\}$
		and $\bdry_Z W = \closure(X-X_\geq) \times \{1\} \cup \bdry X \times I$ and $F = f \times \id_I$.
	\end{Proof}
	
	In summary we have obtained the following Theorem:
	
	\begin{Theorem}\label{thm:Bordism}
		Given a compatible or auxiliary condition $\mathcal{C}$, then
		the (oriented) bordism-groups of MHSS satisfying $\mathcal{C}$ are well-defined.
		Further,
		\begin{enumerate}[(i)]
			\item The graded abelian group $\Omega_*^{\mathcal{C}}(Z,W)$
			is a graded module over (oriented) manifold-bordism coefficients
			$\Omega_*^{\txt{M(S)TOP}}(\{\txt{pt}\})$.
			\item Given an (oriented) $X$ with $\mathcal{C}(X)$ and $\mathcal{C}(\bdry X)$,
			there is a geometric fundamental-class
			given by $[X]:=[X \xrightarrow{\id} X]\in\Omega^{\mathcal{C}}_n(X,\bdry X)$, which is
			invariant under (orientation-preserving) stratified homeomorphisms
			$h: (X,\bdry X) \rightarrow (Y, \bdry Y)$
			so $h_*([X]) = [Y]$.
			\item If closed (oriented) spaces $X$ with $\mathcal{C}(X)$ have an invariant
			$\sigma(X)\in G$ in some abelian group $G$,
			such that $\sigma(\bdry X) = 0$ is trivial on boundaries
			and $\sigma(X\sqcup (-Y)) = \sigma(X) - \sigma(Y)$,
			then there are induced group-homomorphisms
			$\hat{\sigma}_n : \Omega^{\mathcal{C}}_n( Z ) \rightarrow G$
			given on fundamental-classes by $\hat{\sigma}([X]) = \sigma(X)$.
			\item If $\mathcal{C}$ is compatible, then this bordism-theory is a
			generalized homology-theory.
		\end{enumerate}	
	\end{Theorem}
	
	\begin{Proof}
		For part (iii), simply define $\hat{\sigma}_n( [f: X \rightarrow Z] ) := \sigma(X)$.
		For the definition of these bordism-groups see Lemma \ref{lemma:BordismDef} and for (i)
		Lemma \ref{lemma:BordismProperties}, which together with Lemma \ref{lemma:Excision} also
		shows (iv). Part (ii) follows from the definition of the bordism-groups and by
		$\cyl(h)$ (with the collapse-map to $Y$) 
		being a bordism from $h_*([X])$ to $[Y]$, because
		it is (orientation-preserving) stratified homeomorphic to $X \times I$ with
		$\mathcal{C}(X) \Rightarrow \mathcal{C}(X \times I)$.
		
		\emph{Remark:}
		This last statement is actually non-trivial, and holds here only, because -- besides orientation -- no additional
		structure is introduced on spaces, so $h$ need not preserve anything
		beyond topology and orientation. Compare this for example to the pl-case, where
		one would have to require $h$ be a pl-homeomorphism to make this statement.
	\end{Proof}
	
	\begin{Example}
		By Lemma \ref{lemma:ExamplesOfCompatibleConditions}\,ii,
		one could fix a maximal number of meeting strata (the total number of strata is not local, and indeed not preserved under
		gluing), in potentially fixed codimensions and combine this auxiliary condition with any compatible condition
		into a compatible condition.
		
		The resulting bordism-theories have canonical maps to the theory obtained for the original compatible condition,
		but will not typically have (unstratified) homeomorphism-invariant fundamental-classes:
		E.\,g.\ start with a (non null-bordant) manifold $M$ and stratify it as $X = M \supset B$ with a
		(non null-bordant) submanifold $B$ (e.\,g.\ $\mathbb{CP}^2 \subset \mathbb{CP}^4$),
		then bordisms $M \bordant X$ will typically have to \inquotes{terminate} the stratum $B$
		with a new stratum of dimension $\dim(B) - 1$ (cf.\ the constructions of §\ref{sec:IntrinsicTransversality}),
		which is not allowed here.
	\end{Example}

	\section{Sufficient Condition on Links}\label{sec:AlgTopCond}
	
	This section is primarily concerned with proving the following theorem, giving a sufficient
	-- and pointwise, thus auxiliary (cf.\ Lemma \ref{lemma:ExamplesOfCompatibleConditions})
	-- condition for spaces satisfying the
	requirements of the geometric theory of this paper.
	Similar algebraic conditions occur for example in \citep{QuinnEndsOfMaps1,RanickiYamasakiControlledL}.
	
	\begin{Theorem1}
		Given a $5$-gapped MHSS $(X,\bdry X)$ with simple links,
		then $X$ has sufficiently many relative cuts.
	\end{Theorem1}
	
	\begin{Corollary}\label{cor:AlgTopCond}
		In the language of §\ref{sec:CompatibleConditions}:
		Being a $5$-gapped MHSS with simple links is compatible.
	\end{Corollary}
	
	\begin{ProofOfCorollary}
		The hypotheses depend only on codimensions and holink-fibers, thus are auxiliary by
		Lemma \ref{lemma:IntercompatibilityOfConditions} and therefore any cut is automatically
		a $\mathcal{C}$-cut.
		Cuts always exist, by applying the theorem inductively over skeleta, in each step, first to the boundary,
		extending as a product into the boundary-collar and then extending to the rest of $X$.
	\end{ProofOfCorollary}
	
	\begin{Remark}\label{rmk:Codim4}
		For pairs with $\pi$ additionally \inquotes{good} (e.\,g.\ poly-(finite\Slash{}cyclic) \citep{FreedmanQuinnTop4}),
		the gap-hypothesis can be weakened to a difference of $4$, see \citep{QuinnEndsOfMaps3}.
		In this case the 4-skeletons need no longer be locally-conelike however,
		which might lead to issues when applying the material in
		§\ref{sec:IntrinsicTransversality} and §\ref{sec:LClasses}.
	\end{Remark}

	\begin{Definition}
		A \defName{mapping-cylinder neighborhood}, is an end-completion in the sense of
		\citep{QuinnEndsOfMaps1,QuinnEndsOfMaps2}, this means:
		
		First, the upper stratum $M$ of $X = M \cup B$ can be completed
		to a compact manifold-with-boundary $M'$ with $M'-M \subset \bdry M'$.
		I.\,e.\ we can essentially add $\bdry_0 M' = \closure(M' - M)$ to the boundary,
		making the upper stratum compact.
		
		Next, we require the existence of a neighborhood $U$
		of the lower skeleton $B$,
		and a continuous map
		$g : \bdry_0 M' \rightarrow B$,
		such that $U$ is homeomorphic to $\opencyl( g )$ strictly
		rel $B$ (i.\,e.\ the homeomorphism from $U$ to $\opencyl( g )$
		maps points in $B$ via $\id_B$, while points of $U-B$
		are mapped to the cylinder-part $\bdry_0 M' \times (0,1) = \opencyl( g ) - B$).
		
		The neighborhood $U$ is than the mapping-cylinder neighborhood of $B$.
	\end{Definition}
	
	Much of the difficulty is in obtaining the \emph{manifold} $\bdry M'$.
	Based on controlled topology results of Quinn \citep{QuinnEndsOfMaps2} one obtains:
	
	\begin{Lemma}\label{lemma:ObstructionGroupsVanish}
		Given a MHSS $X$ with empty boundary $\bdry X=\emptyset$ and with simple links,
		then the obstruction-groups of \citep{QuinnEndsOfMaps2}, $H_0^{\txt{lf}}(X^k; \mathcal{S}(p_k))$
		of obstructions to the existence of cylinder-neighborhoods of $X^k \subset X^{k+1}$ (controlled end-completions) and
		$H_1^{\txt{lf}}(X^k; \mathcal{S}(p_k))$ to their uniqueness (controlled h-co\-bord\-isms) vanish for all $k$.
	\end{Lemma}
	\begin{Proof}
		The vanishing of end-obstruction groups in this case can be deduced either
		from the Atiyah--Hirzebruch like spectral sequence
		\citep[Thm.\,8.7 (p.\,423)]{QuinnEndsOfMaps2} or directly from the properties of the obstruction-spectrum (making it a
		\inquotes{homology}-theory in the sense of \citep[Thm.\,8.5 (p.\,421)]{QuinnEndsOfMaps2})
		shown in \citep[Thm.\,5.11 (p.\,400)]{QuinnEndsOfMaps2}:
		Inductively over skeleta, use that skeleta are \inquotes{pure},
		thus $p$-NDR (co-fibrations) so that restriction induces a homotopy-fibration,
		together with the result that over individual strata, the $\Wh(\pi \times \Integers^k) = 0$ condition
		ensures vanishing of the obstruction groups
		\citep[Thm.\,2.1 (p.\,282)]{QuinnEndsOfMaps1} (applied with $D = \emptyset$, $\bdry X = \emptyset$).
	\end{Proof}
	
	\begin{Lemma}\label{lemma:CylExists}
		Let $X$, $\bdry X=\emptyset$, be a MHSS with simple links.	
		Then given $k>0$, the $k$-skeleton $X^k$ has a mapping-cylinder neighborhood
		in the next larger stratum $X_l$, if that stratum has dimension $l \geq 6$.
	\end{Lemma}
	
	\begin{Proof}	
		Since the obstruction-groups vanish (Lemma \ref{lemma:ObstructionGroupsVanish}), and $l\geq 6$, the end-theorem,
		formulated most conveniently for the present situation
		in \citep[Thm.\,1.7 (p.\,446)]{QuinnHSS}, applies to yield a mapping-cylinder neighborhood.
	\end{Proof}
	
	\begin{Remark}\label{rmk:Choices}
		The choice of this mapping-cylinder neighborhood is unique only up to controlled h-cobordism
		\citep[Thm.\,1.1b (p.\,357))]{QuinnEndsOfMaps2},
		which under the hypothesis used here are classified by the trivial group. However in general,
		for example when using pl-structures, there is technically a non-trivial choice implicit here.
		
		For an example where the choice of end-completion changes the pl-homeo\-morph\-ism-type,
		see \citep{MilnorHomeoNotPL}. So this is not in contradiction to the pl-results of \citep{BanaglLauresMcClure}.
	\end{Remark}
	
	These mapping-cylinders were understood to be particularly nice (the generating map has an approximate lifting property)
	by Hughes \citep{HughesCylinders}, which was then employed by Connolly and Vajiac
	\citep{ConnollyVajiac} to give a suitable
	homotopy-transversality statement (they also understood
	and used the connection to Quinn's controlled topology-results,
	however with different goals, which makes their algebraic hypothesis unsuitable for the present treatment):
	
	\begin{Lemma}\label{lemma:HoCutExists}
		Let $X$ be a MHSS with empty boundary.
		Then given $k>0$, such that the $k$-skeleton $X^k$ has a mapping-cylinder
		neighborhood in the next larger stratum $X_l$,
		and a cut of the $k$-skeleton, then this cut can be extended as a homotopy-cut
		(see proof and next lemma) into this mapping-cylinder-neighborhood.
	\end{Lemma}
	
	\begin{Proof}
		A homotopy-cut, here means, a stratified (codimension $1$) MHSS subspace $X_{\{0\}}$, which is a stratified neighborhood-deformation-retract.
		This result can be found in \citep[Thm.\,2.2 (p.\,529)]{ConnollyVajiac}.
	\end{Proof}
	
	Again using controlled topology \citep{QuinnEndsOfMaps2}, this can be seen to automatically induce (in high dimensions)
	a geometric cut (i.\,e.\ a cut in the sense of Def.\,\ref{def:Cut}), this was also noted
	by \citep{ConnollyVajiac}.
	
	\begin{Lemma}\label{lemma:HoCutGeom}
		Let $X$ be a MHSS with empty boundary.
		Then given $k>0$, and a cut of the $k$-skeleton $X^k$, which extends as a homotopy-cut into an open neighborhood $U$
		of $X^k$ in the next larger skeleton $X^l$,
		then, if $l \geq 5$, a (geometric) cut extending the one in the $k$-skeleton,
		including boundary-collars (as germs, i.\,e.\ potentially thinner) into $U$ exists.
	\end{Lemma}
	
	\begin{Proof}
		Note, that while \citep[Thm.\,1.2 (p.\,444)]{QuinnHSS} gives only the absolute case,
		the proof is actually given inductively via the relative result \citep[First Lemma of §5.3 (p.\,492)]{QuinnHSS},
		which is the relative result stated above.
	\end{Proof}
	
	Finally we can assemble these ingredients to prove Theorem \ref{thm:AlgTopCond}.
	
	\begin{ProofOfThm1}
		We first treat the case $D = X^k$ and $\bdry X = \emptyset$.
		Given a cut between $D \cap B_\pm$,
		we show a cut of $X^{k+1}$ rel $X^k$ exists (between $X^{k+1} \cap B_\pm$),
		rel $D$. This then proves the theorem, by inductive application to skeleta of $X$.	
		Let $S$ be a connected component (meeting $X^k$) of the next larger stratum $X_{k+1}$.
		
		Case $\dim(S) \leq 6$: Since $X^0$ consists of points, these can be absorbed into the closed disjoint sets $B_\pm$
		retaining them being closed disjoint. Thus w.\,l.\,o.\,g.\ we may assume $X^0 = \emptyset$.
		Hence, by the gap-hypothesis, $S$ is a manifold, thus can be cut using manifold-transversality.
		
		Case $\dim(S) \geq 6$: Using the controlled topology (end-theorem) results of Quinn \citep{QuinnEndsOfMaps2}
		via Lemma \ref{lemma:CylExists},
		a mapping-cylinder-neighborhood of $X^k \subset S \cup X^k$ exists.
		By the construction of \citep{ConnollyVajiac} via the cylinder-theorem of \citep{HughesCylinders}
		in Lemma \ref{lemma:HoCutExists},
		the cut of $X^k$ extends as a homotopy-cut into this cylinder-neighborhood $U$ of $X^k\subset S \cup X^k$.
		Again using a controlled topology result (h-cobordism) through Lemma \ref{lemma:HoCutGeom},
		we find a geometric cut extending the one of $X^k$ into $U \subset S \cup X^k$.
		For the extension to the remainder of $X^{k+1}$, see the final step of the relative case below.
		
		Next we treat the case of $D = X^k \cup C$ for an open subset $C \subset X$, for $\bdry X = \emptyset$.
		We are given	an open neighborhood $U$ of $X-C$.
		
		Case $\dim(S) \leq 6$: The argument above still applies, using relative mani\-fold-transversality.
		
		Case $\dim(S) \geq 6$: 
		We are given a cut $(C_\geq, C_\leq,C_{\{0\}})$ on the open $C$ and an open neighborhood $U$ of $X-C$.
		We have to extend this cut to $X$ while leaving it unaltered on $X-U \subset C$.
		To this end, we will be \inquotes{gluing} the new part on $U$ along $U\cap C$ to the old part on $C$.
		We can shrink $U$ to a smaller neighborhood $U' \subset U$ of $X-C$,
		then the statement for $U'$ implies the one for $U$.
		
		Also note, that we can essentially focus on $C_{\{0\}}$: Since $\bdry X= \emptyset$,
		we have $\bdry C_\geq = \bdry C_\leq = C_{\{0\}}$ and $\bdry C_{\{0\}} = \emptyset$,
		so boundary collars come from a bi-collar of $C_{\{0\}}$.
		Further, potentially by making the neighborhood $U$ of $X-C$ smaller (see above),
		we can assume $(C_{\{0\}})^0 = \emptyset$ as
		we have to produce some open overlap with the relative part over $C$ to glue the new part of the cut
		to the old part, but the discrete point-set $(C_{\{0\}})^0$ can be absorbed either in the part that is identically
		kept as is, or ignored where the new part is created.
		Hence if $\dim(S) = 6$ the gluing to the relative part is again the manifold-case as treated above.
		
		We can thus assume $\dim(S) \geq 7$ and therefore $\dim(S \cap C_{\{0\}}) = \dim(S)-1 \geq 6$.
		Further $C_{\{0\}}$ again has simple links, as this is an auxiliary hypothesis.
		In particular the skeleton $C_{\{0\}}^{k-1} = X^k \cap C_{\{0\}}$ has a cylinder-neighborhood in the
		stratum $S \cap C_{\{0\}}$ corresponding to $S$, say $\opencyl( f_0: \bdry M^C_0 \rightarrow C_{\{0\}}^{k-1} )$.
		Thus, close to $X^k$ the bi-collar of $C_{\{0\}}$ is of the form
		$\opencyl( f_0 ) \times \Reals = \opencyl( f_0 \times \id_{\Reals} )$.
		
		Next, we produce an \inquotes{absolute} cut $X'_{\{0\}}$ (by the first part of the proof) of $X^k \cup S$,
		extending the one of $X^k$.
		Using the same logic as for $C_{\{0\}}$, we find that near this cut $X^k \cup S$
		is of the form $\opencyl( g_0 \times \id_{\Reals} )$.
		By \citep[Thm.\,1.1\,b (p.\,357)]{QuinnEndsOfMaps2}
		the completion over the overlap $U \cap C$ is unique up to h-cobordism.
		As explained in more detail in \citep[p.\,1149]{QuinnEndsOfMaps4},
		this means, for two completion $V = \opencyl(f) = \opencyl(g)$, fixing $\delta > 0$, such that for
		$t,t'$ chosen small enough, the region of $V$ between $\{ [x,\tau] \in \opencyl(f) | \tau = t \}$
		and $\{ [x,\tau] \in \opencyl(g) | \tau = t' \}$ is a ($\delta$,h)-cobordism.
		
		By \citep[Thm.\,1.1\,c (p.\,357)]{QuinnEndsOfMaps2}, the obstruction for this ($\delta$,h)-cobordism
		to being an $\epsilon$-product (we technically first fix an $\epsilon >0$, then the approximate h-cobordism-theorem
		\citep[Thm.\,1.6 (p.\,359)]{QuinnEndsOfMaps2} yields the appropriate $\delta$ to use in the step above),
		is in $H_1^{\txt{lf}}(U \cap C; \mathcal{S}(p_k|_{U \cap C})) = 0$, which vanishes by Lemma
		\ref{lemma:ObstructionGroupsVanish} (the open subspace of $X^k \cup S$ over $U\cap C$ has simple-links,
		because $X$ does).
		See Remark \ref{rmk:CutUniqueness} below for an explanation, why this works,
		while extension of collars may not work.
		
		The technical approach we take to fit the cuts together is similar in a sense to the one of \citep{ConnollyVajiac}:
		First note, if we write coordinates in the cylinders as
		$([x,t],s)\in \opencyl(f_0 \times \id_{\Reals}) = \opencyl( f_0 ) \times \Reals$,
		then $C_{\{0\}}$ consists precisely of the points with $s=0$, but this can isotoped by $\alpha_t$ (see below)
		rel $X^k$ ambient in the part $M$ of $S$ over $C$, to
		\begin{equation*}
			C^\chi_{\{0\}} = \big\{ ([x,t],s)\in \opencyl( f_0 ) \times \Reals \big| t=\chi(f_0(x),s) \big\}
		\end{equation*}
		for $\chi \in (0,1)$ (the choice of $\chi$ will be explained in the next paragraph).
		The isotopy $\alpha_t$ can be constructed
		(at least until $\chi =1$, but parts far away from $X^k$ can later be repaired using manifold-transversality on $S$ instead,
		see below) as follows:
		We may assume $\chi(y,s) < s$, because $\chi$ needs only to be small enough (see next paragraph).
		To obtain $\alpha$, first push along the cylinder-coordinate \inquotes{t} up to $t=s$, thereby moving
		$C^\chi_{\{0\}}$ to $C^{(y,s)\mapsto s}_{\{0\}}$ rel $X^k$ (because $\chi>0$ and $\chi \rightarrow 0$ for
		$s\rightarrow 0$). This $C^{(y,s)\mapsto s}_{\{0\}}$ 
		(that is, points with $t=s$) can in turn be pushed along the product-coordinate \inquotes{s}
		to $s=0$, i.\,e.\ to $C_{\{0\}}$.	 
		Similarly we get $X^\eta_{\{0\}}$ isotopic to $X_{\{0\}}$ via $\beta_t$.
		
		We pick $\chi$ and $\eta$ as follows:
		We focus on the region over one side of the bi-collar in $X^k$ (that over $s>0$) in the coordinates above and
		let $\epsilon : ( X^k_{\{0\}} \times (0,1) ) \cap (C\cap U) \rightarrow (0,1), (x,s) \mapsto s$.
		This $\epsilon$ can be chosen as a map (rather than a number) to obtain a map $\delta$,
		by repeated application of the relative
		statement of \citep[Thm.\,1.6 (p.\,359)]{QuinnEndsOfMaps2}\footnote{
			Our spaces are locally compact and paracompact, so for every point $x$, there is a compact neighborhood $K_x$
			and for the number $\epsilon_x = \min_{K_x}(\epsilon)$, there is a number $\delta_x$,
			use a decomposition of one to glue these $\delta_x$ into $\delta$.
			Alternatively, use
			\citep[Thm.\,1.1 and 1.2 (p.\,357)]{QuinnEndsOfMaps2} directly,
			see also \citep[\inquotes{Proposition} (p.\,492)]{QuinnHSS}.}.
		This determines (by \citep[Thm.\,1.6 (p.\,359)]{QuinnEndsOfMaps2}) a map
		$\delta : ( X^k_{\{0\}} \times (0,1) ) \cap (C\cap U) \rightarrow (0,\infty)$,
		which in turn determines $\chi, \eta$, small enough,
		that the region between
		$\{ [x,\tau] \in \opencyl(f_0 \times \id_\Reals) | \tau = \chi(f_0\times \id(x)), f_0\times\id(x)\in X^k_{\{0\}} \times (0,1) \}$
		and
		$\{ [x,\tau] \in \opencyl(g_0 \times \id_\Reals)|\tau = \eta(f_0\times \id(x)), g_0\times\id(x)\in X^k_{\{0\}} \times (0,1) \}$
		is a ($\delta$,h)-cobordism
		(again, as for obtaining $\delta$ from $\epsilon$ as map, we obtain $\chi$ and $\eta$ from $\delta$,
		corresponding to the numbers $t$ and $t'$ above).
		This is an $\epsilon$-product by choice of $\delta$ and vanishing of the obstruction-group,
		see above.
		We define an isotopy $\gamma_t$ of $M$ by pushing along this product coordinate, moving $X_{\{0\}}^\eta$ to
		$C_{\{0\}}^\chi$ (as a \inquotes{buffer} to start\Slash{}stop the push, use that the cylinder-coordinates
		of $\opencyl(f_0\times \id_\Reals)$ and $\opencyl(f_0\times \id_\Reals)$
		turns the regions between $\onehalf \eta$ and $\eta$ and $\chi$ and $2\chi$ respectively into products).
		
		By choice of $\epsilon = s$, the resulting $\epsilon$-product (and hence the isotopy)
		gets arbitrary small (after projection with the control-map, which is simply the
		cylinder-collapse here) when approaching the middle ($s=0$)
		of the cut. Thereby, when approaching the lower stratum (which happens by the choices of
		$\eta$ and $\chi$ for $s \rightarrow 0$), this isotopy
		moves the cut in the upper stratum in way, such that the 
		center (later to become $X^{k+1}_{\{0\}}$) always (i.\,e.\ for all $t\in[0,1]$) fits together with the center
		$X^k_{\{0\}}$ of the cut in the lower stratum.
		
		To extend this away from $C$, use a map $d: X^k \rightarrow [0,1]$ such that $d=1$ on $X-U$ and $d=0$ on $X-C$
		to multiply the parameter on the isotopy above to get an isotopy on $C$,
		that agrees with the one above on $X-U$, and with the identity near $X-C$. Extend this isotopy by the identity
		on $X-C$. We call the extended result $\bar{\gamma}_t := \gamma_{d t}$.
		Similarly $\alpha_{d t}$ can be extended by the identity to the remainder of the cylinder-neighborhood.
		Call this extended isotopy $\bar{\alpha}_t := \alpha_{d t}$
		We thus obtain $X_{\{0\}} := ( \bar{\alpha}^{-1} * \bar{\gamma} * \beta)_1 (X'_{\{0\}})$
		as the image	of $X'_{\{0\}}$ under the composed isotopy at $t=1$.
		
		This is stratified \BiCollared{} rel the \BiCollar{} given in $X^k$, because $X'_{\{0\}}$ is and the isotopy is rel $X^k$.
		Further, on $C-U$, it agrees with $C_{\{0\}}$, by construction of the isotopy.
		
		Finally, apply relative manifold-(map-)transversality to extend this cut to all of $S \cup X^k$:
		Use (the inverse of) the bi-collar $c: X_{\{0\}} \times (-\epsilon,\epsilon) \rightarrow U$ to
		define $g: \im(c) \rightarrow (-\epsilon,\epsilon)$ as the projection to the bi-collar coordinate
		and extend this $g$ continuously to $S \cup X^k$ such that $g|_{B_\pm} = \pm 1$
		and $g^{-1}(\{0\}) \cap U = c(X_{\{0\}}\times\{0\})$ by regularity of the metric $X$.
		Since $X$ is metric, thus normal, there is an open neighborhood $V$ of $X^k$ with $\closure(V) \subset U$,
		make $g|_S$ transverse to $0$ rel $g^{-1}(\{\pm 1\}) \cup V$ \citep{KirbySiebenmannEssays,QuinnTransversality}.
		To this end, note, that $g|_S$ is already transverse to $0$
		on $U$ by construction. This yields $g^S_\perp$, which fits together with the original $g$ on $V$ to form $g_\perp$
		on $X^k \cup S$.
		The new cut of $X^k \cup S$ is $(X_\geq,X_\leq, X_{\{0\}}) := g_\perp^{-1}( [0,1], [-1,0], \{0\} )$.
		Note, that $X$ does not have a boundary,
		and hence $\bdry X_\geq = \bdry X_\leq = X_{\{0\}}$ so that the bi-collar of $X_{\{0\}}$ obtained by the
		transversality-arguments above provides boundary-collars for $X_\geq$ and $X_\leq$, and this is indeed a cut
		in the sense of Def.\,\ref{def:Cut}.
		Different components $S'$ of $X_{k+1}$ can be treated independently, so we may assume
		this is a cut of $X^{k+1}$ (rel $X^k$).
		
		Case $\bdry X \neq \emptyset$:
		The hypothesis of the theorem is auxiliary, hence satisfied by $\bdry X$.
		Apply the case for empty boundary above to $\bdry X$,
		to obtain a cut of the boundary between $\bdry X \cap B_\pm$ and rel $D \cap \bdry X$.
		Extend this as a product to the boundary-collar of $\bdry X \subset X$,
		choosing the collar \inquotes{thin} enough, so this is still between $B_\pm$
		and a product with the collar on the boundary of $X^k$ (the relative cut was assumed boundary-compatible).
		Next apply the case for empty boundary to $X-\bdry X$ rel $(D \cup \txt{boundary-collar})$
		and for a $U$ with $U \cap \bdry X = \emptyset$.
	\end{ProofOfThm1}
	
	\begin{Remark}\label{rmk:CutUniqueness}
		The proof above only uses, that the cylinder-neighborhoods $\opencyl(f_0 \times \id_{\Reals})$ and
		$\opencyl(g_0 \times \id_{\Reals})$ can be isotoped (in $S$ rel $X^k$) to match a single \inquotes{level} $t$.
		Matching all levels at once, would require the control-map to the holink (see \citep{QuinnHSS})
		to be a $(\delta,2)$-equivalence (it generally only is a $(\delta,1)$-equivalence) 
		or holink-fibers would need to be simply-connected, see
		\citep[Prop.\,2.4 (p.\,1150) and Cor.\,1.2 (p.\,1143)]{QuinnEndsOfMaps4}.
		However, here we explicitly ask only for the cut and \emph{some} collar,
		not necessarily an extension of the collar of $C_{\{0\}}$.
		This is essentially the difference between h-cobordism triviality
		as opposed to its uniqueness up to pseudo-isotopy, which is what requires the strengthened hypothesis
		in \citep{QuinnEndsOfMaps4}.
		
		The absolute version of the theorem is essentially a map-transversality theorem to a target with trivial normal-bundle.
		A fully (including collars) relative statement would immediately yield a map-transversality theorem
		to arbitrary normal-bundles. This seems to be true for simply-connected holink-fibers by
		\citep[Cor.\,1.2 (p.\,1143)]{QuinnEndsOfMaps4}, but interestingly may not hold for all \inquotes{simple links}
		by the issues pointed out above.
	\end{Remark}
	
	\begin{Remark}\label{rmk:TransversalityAndEnds}
		There is an intimate relation of end-completion, h-co\-bord\-isms and transversality
		\citep{ConnollyVajiac}\citep[p.\,14f]{QuinnEndsOfMaps2} in general.	
		Interestingly, when posed in the context of an inverse-excision\,/\,cutting-problem,
		here stratified-trans\-vers\-ality is slightly \inquotes{easier} in principle, as there is an additional degree of freedom in slightly "moving" the cut around
		(its position is essentially fixed \inquotes{homotopically only}).
		This allows for example the removal of \inquotes{strange points}.
		This \emph{might} also make the dimension (gap) requirements less severe,
		but to establish a formal result on this seems technically extremely demanding.
	\end{Remark}

	\section{Simultaneous Transversality}\label{sec:IntrinsicTransversality}
	
	Quinn \citep{QuinnIntrinsicSkeleta} shows -- for the purpose of showing that
	intersection homology is topologically invariant on MHSS --
	that there are certain intrinsic skeleta of MHSS, which we will denote by $|X|$.
	In \citep{QuinnIntrinsicSkeleta}
	these are denoted $X_{0,0}$ (cf.\ Rmk.\,\ref{rmk:IntrinsicConstruction}),
	while $|X|$ denotes the topologically intrinsic stratification, however in the cases
	in which we are interested both agree, see (Q2) below.
	These depend only on the (unstratified) homeomorphism-type of a space $X$.
	Generally $|X|$ is only a HSS (its strata may not be manifolds) however
	by \citep[Thm.\,1.1 (p.\,235)]{QuinnIntrinsicSkeleta} the following properties are satisfied,
	where (Q2) also clarifies the manifold-stratum problem:
	\begin{enumerate}[(Q1)]
		\item If $Y$ is another stratification as MHSS of the underlying topological space of $X$, then
		$|X|$ is a coarsening\footnote{I.\,e.\ the identity $\id : Y \rightarrow |X|$ is a coarsening, meaning
			preimages of components of strata of $|X|$ are unions of components of strata of $Y$.
			In other word, $Y$ is a \inquotes{subdivision} of $|X|$.} of $Y$.
		In particular, if $X$ and $Y$ are (unstratified) homeomorphic, then $|X|$ and $|Y|$ are stratified homeomorphic
		(the homeomorphism and its inverse are both coarsenings).
		\item If either $X^4$ is locally conelike or $k\geq 5$, then the $k$-stratum of $|X|$ is a manifold and the $k$-skeleton
		agrees (potentially after re-indexing) with the one of the topologically intrinsic stratification
		given by equivalence classes of points with homeomorphic neighborhoods.
		\item If $|X|^0 = \emptyset$, then $|X \times \Reals| = |X| \times \Reals$.
	\end{enumerate}
	
	\begin{Remark}\label{rmk:IntrinsicConstruction}
		The construction of $|X| = X_{0,0}$ in \citep{QuinnIntrinsicSkeleta} is done roughly as follows:
		By a double induction (over $j$, then $k$ see below) from top to bottom
		over strata, components of the $k$-stratum are \inquotes{promoted} to the $j$ stratum
		whenever the link has the correct homotopy-type
		for the result to be a manifold and the result would still be an HSS
		(plus some technical tweaks to make things work on a formal level).
		
		It is sometimes helpful for the understanding of some of the constructions, to be aware,
		that $|X|$ is obtained by
		promoting \emph{components of strata} in $X$ to become parts of higher strata.
	\end{Remark}
	
	\begin{Remark}\label{rmk:StrangePoints}
		By property (Q3), it is already evident, that this is indeed a problem near \emph{points} only.
		By the discussion of \citep[Thm.\,1.1 (p.\,235)]{QuinnIntrinsicSkeleta}, these problems arise at and only at
		points that are not \inquotes{1-LC} (the complement has non-trivial local fundamental-group near the end),
		while being promotable (Rmk.\,\ref{rmk:IntrinsicConstruction}) into higher strata, which requires
		in particular, that holink-fibers are homology-spheres see \citep{QuinnIntrinsicSkeleta}.
		
		More concretely, this issue occurs for example for cones on non-simply-connected homology-spheres,
		which cannot be manifolds near the cone-point, but whose product with the real line is a manifold
		(essentially by Edwards' double-suspension-theorem, see \citep{EdwardsDoubleSuspension}),
		so the cone-point must be singular in $|X| \times \Reals$ (for the locally conelike $|X|$ to have manifold strata),
		but it must be non-singular in $|X \times \Reals|$ (for this to have the topologically intrinsic stratification).
	\end{Remark}
	
	By (Q3), the intrinsic space of a boundary-collar is a collar of the intrinsic space of the boundary,
	except possibly at $0$-strata, that come from cutting $1$-strata, where that $1$-stratum is not a genuine singularity --
	it gets absorbed as part of a higher-dimensional stratum of the intrinsic space -- but the $0$-stratum is.
	The central idea of this section is, that since this is not a problem after stabilization by $\Reals$ (see lemma below),
	as long as both $X$ and $|X|$ have sufficiently many relative cuts, the cylinder in between can be stratified and
	cut in a compatible way.
	
	\begin{Lemma}\label{lemma:SimCut}
		Given a MHSS $X$ with locally conelike $3$-skeleton and
		a bound\-ary-compatible cut $(X_\geq,X_\leq,X_{\{0\}})$ of $X$,
		the triple $(|X_\geq \times \Reals|,|X_\leq \times \Reals|,|X_{\{0\}} \times \Reals|)$
		is a boundary-compatible cut of $|X \times \Reals|$.
	\end{Lemma}
	
	\begin{Proof}
		Since $X$ was assumed to have locally conelike $3$-skeleton,
		$X \times \Reals$ has locally conelike $4$-skeleton and thus $|X \times \Reals|$
		has manifold-strata.
		
		$(\bdry X_\geq \times \Reals)^0 = \emptyset$, thus given a boundary-collar
		$\bdry X_\geq \times [0,\infty) \subset X_\geq$ the corresponding
		$|\bdry X_\geq \times [0,\infty) \times \Reals| = |\bdry X_\geq \times \Reals| \times [0,\infty)$,
		using (Q3) and a double-space $X_\geq \cup_{\bdry X_\geq} X_\geq$,
		is a boundary-collar of $|\bdry X_\geq \times \Reals| \subset |X_\geq \times \Reals|$.
		Similarly for $X_\leq$.
	\end{Proof}
	
	\begin{Definition}\label{def:CoarseGluing}
		Given MHSS $X^{n+1}, Y^{n+1}, Z^n$, with $Z \subset \bdry X$
		and $X, Z$ compact,
		where $Y$ has an end stratified homeomorphic to $|Z \times \Reals|$
		(i.\,e.\ there is an open stratified embedding $|Z \times (0,\infty)| \subset Y$ with
		the restriction to $|Z \times (0,1]|$ proper\footnote{Compact subsets of $Y$ restricted to
			$|Z \times (0,1]|$ are thus compact, so $Y$ is open towards $0$ and effectively \inquotes{missing a boundary}.},
		and with $Y-|Z \times \Reals|$ compact).
		Then the \defName{intrinsic gluing $X \GlueIntrin{Z} Y$} along $Z$ is the MHSS (see
		Lemma \ref{lemma:HalfIntrinsicCyl} showing this is indeed a MHSS)
		given by stratifying the adjunction-space
		$X \cup_Z (Y \cup (Z \times \{0\}))$ as:
		\begin{equation*}
			W^k := X^k \cup Y^k \cup
			\begin{cases}
				\emptyset & \txt{if }k\geq n	\\
				Z^k	& \txt{otherwise}
			\end{cases}
		\end{equation*}
		
		We call a cut of such an intrinsic-gluing \defName{middle-compatible},
		if it is bound\-ary-compatible, restricts to a bound\-ary-compatible cut of $X$ and
		near $Z$ is given as in Lemma \ref{lemma:SimCut} on non-top strata of $Z$
		(i.\,e.\ it is essentially \inquotes{boundary-compatible}
		when restricted to $Y$, except that the boundary at $Z$ is missing).
	\end{Definition}
	
	\begin{Lemma}\label{lemma:CutIntrinsicGluings}
		Given $5$-gapped MHSS $X,Y, Z$ with simple links,
		and such that the intrinsic gluing (Def.\,\ref{def:CoarseGluing}) $X \GlueIntrin{Z} Y$ is defined,
		then $X \GlueIntrin{Z} Y$ has sufficiently many absolute cuts, such that the resulting spaces and their boundaries
		are again of this form.
		The cut can be chosen middle-compatible and rel a middle-compatible cut of $\bdry ( X \GlueIntrin{Z} Y )$.
	\end{Lemma}
	
	\begin{Proof}
		Let $B_\pm \subset X \GlueIntrin{Z} Y$ closed, disjoint.
		
		We cut $Z$
		between $Z\cap B_\pm$ using Theorem \ref{thm:AlgTopCond} extending the given cut on
		$\bdry Z \subset \bdry ( X \GlueIntrin{Z} Y )$.
		
		$Z$ has, by Def.\,\ref{def:CoarseGluing}, a neighborhood in $X \GlueIntrin{Z} Y$
		that intersects $Y$ in $|Z \times (0,\infty)|$ in $Y$,
		by Lemma \ref{lemma:SimCut}, $|(Z_\geq,Z_\leq,Z_{\{0\}})\times \Reals|$
		is also a cut of this neighborhood $|Z \times \Reals|$.
		The relative cut is of this form, if we are not in the top-stratum, by hypothesis (it is middle-compatible).
		
		If we are not in the top stratum,
		by Theorem \ref{thm:AlgTopCond}, the cut of $\bdry Y$ near $\bdry Z$ can be extended to one of $\bdry Y$,
		then to one of $Y$ rel boundary.
		Similarly, on $X$ extend the cut using Theorem \ref{thm:AlgTopCond}.
		
		If we are in the top-stratum,
		we still obtain a cut like this, but it may not agree with the one of the boundary near the top-stratum of $Z$.
		However, since the cut is boundary-compatible, we can embed a product $\bdry ( X \GlueIntrin{Z} Y ) \times [0,1]$
		near $\bdry ( X \GlueIntrin{Z} Y )$ via the boundary-collar, s.\,t.\ $\bdry ( X \GlueIntrin{Z} Y ) \times \{0\}$
		is the intersection with the boundary.
		Use the isotopy constructed in the proof of Theorem \ref{thm:AlgTopCond} (to fit together relative cuts),
		for the cuts of $Z \times \Reals$ in the boundary and let the $t$-coordinate of this isotopy run with the
		value in $[0,1]$ of the embedding above to fit both together.
		
		These cuts of $Y$ and $X$ again form a intrinsic gluing:
		The cut $(X \GlueIntrin{Z} Y)_{\{0\}}$ near $Z$ in $Y$ was chosen
		as $|Z_{\{0\}} \times (0,\infty)|$, so the end of $Y \cap (X \GlueIntrin{Z} Y)_{\{0\}}$
		is again of the form $|Z_{\{0\}} \times (0,\infty)|$
		i.\,e.\ $(X \GlueIntrin{Z} Y)_{\{0\}}$ is an intrinsic gluing of $5$-gapped MHSS with simple links
		(because this is an auxiliary condition) and similar for $Z_\geq$ and $Z_\leq$.
		So these again satisfy the hypothesis of the lemma.
	\end{Proof}
	
	Being an intrinsic gluing as in Def.\,\ref{def:CoarseGluing} by itself is \emph{not} a
	compatible condition, as gluing two of them is \emph{not} again an intrinsic gluing.
	Hence we need to \inquotes{stabilize} against multiple gluing:
	
	\begin{Definition}\label{def:MultiGluing}
		A MHSS $X$ is an \defName{intrinsic multi-gluing of $5$-gapped MHSS with simple links},
		if there exist finitely many
		$5$-gapped MHSS $X_j,Y_j,Z_j$ with simple links ($j = 0, \ldots, N$),
		such that
		there are $\bar{X}_1 = (X_0 \GlueIntrin{Z_0} Y_0) \cup_{A_1} (X_1 \GlueIntrin{Z_1} Y_1)$,
		with $A_1 \subset \bdry (X_0 \GlueIntrin{Z_0} Y_0), \bdry (X_1 \GlueIntrin{Z_1} Y_1)$
		a MHSS (automatically $5$-gapped with simple links),
		$\bar{X}_2 = \bar{X}_1 \cup_{A_2} (X_2 \GlueIntrin{Z_2} Y_2)$ etc.\ 
		such that $X = \bar{X}_N$.
		
		Some of the $Y_j$ and $Z_j$ might be empty, so this includes gluings with $5$-gapped MHSS with simple links
		and disjoint unions.
	\end{Definition}
	
	\begin{Lemma}\label{lemma:IntrinsicGluedCompletion}
		Being an intrinsic multi-gluing of $5$-gapped MHSS with simple links
		is a compatible condition in the sense of §\ref{sec:CompatibleConditions}.
	\end{Lemma}
	
	\begin{Proof}	
		Products with manifolds (first for intrinsic gluings):
		$(X \GlueIntrin{Z} Y) \times M = (X\times M) \GlueIntrin{Z\times M} (Y \times M)$
		are again intrinsic gluings: The end of $Y$ being of the form $|Z \times (0,\infty)|$ has empty zero-skeleton,
		and by locality, via (Q2) and the topologically intrinsic stratification, the product with $M$
		is essentially a product with $\Reals^m$ and thus by (Q3)
		$|(Z\times \Reals^m) \times (0,\infty)| = |Z \times (0,\infty)| \times \Reals^m$ so the form of the end
		of $Y\times M$ is correctly $|(Z\times M) \times (0,\infty)|$.
		The space $\bar{X}$ being glued of such spaces is of course also stable under manifold-products.
		
		Constructing cuts, one needs to be a little careful:
		Boundaries of intrinsic gluings are again intrinsic gluings (this must be allowed, so products by manifolds
		with boundary, in particular with the interval $I$, are allowed),
		however, we only required a structure as intrinsic gluing to \emph{exist},
		so even though the gluing interface (of a standard\Slash{}non-intrinsic gluing) has such a structure,
		the two glued spaces need not \inquotes{agree} on the choice of intrinsic-gluing structure of the boundary,
		which is required however, for the middle-compatible condition on the relative boundary-collars.
		What saves the day, is, that they disagree, at most, on the top-stratum:
		Since the spaces $X$ and $Y$ are $5$-gapped, all non-top strata of $Z$ can be detected as those
		strata of $X \GlueIntrin{Z} Y$ that meet a stratum of (exactly) one dimension higher.
		The top-stratum is treated differently in Lemma \ref{lemma:CutIntrinsicGluings}, so this is allowed.
		
		Hence, cuts exist via Lemma \ref{lemma:CutIntrinsicGluings}, which provides (middle-compatible) cuts,
		first of the boundaries
		between the \inquotes{blocks} of intrinsic gluings, then of the intrinsic gluings themselves rel boundary.
		Finite gluings of spaces of this form are also again of this form (which is what kept intrinsic gluings themselves from
		being compatible).
	\end{Proof}
	
	\begin{Lemma}\label{lemma:BordismSpaceIntrinsic}
		The condition of being an intrinsic multi-gluing of $5$ gapped MHSS with simple links
		has the following properties:
		\begin{enumerate}[(i)]
			\item It is compatible.
			\item Given a $5$-gapped MHSS $X$ with simple links, then
			$X$ is an intrinsic multi-gluing of $5$ gapped MHSS with simple links.
			\item Given a $5$-gapped MHSS $X$ with simple links and no strange points,
			such that $|X|$ is also $5$-gapped with simple links, then
			there is a stratification $W$ of the underlying space of $X\times I$
			with $\bdry W = X \sqcup |X|$ (for orientations, see iv) such that
			$W$ is an intrinsic multi-gluing of $5$ gapped MHSS with simple links.
			\item An intrinsic multi-gluing of $5$-gapped MHSS with simple links $X$,
			has locally conelike $4$-skeleton.
			\item If $X$ has dense-top stratum and no codimension
			$1$ stratum (then so does $|X|$ by (Q1)),
			then $W$ satisfies both as well, in this case, there is a canonical choice for an orientation of $W$, such that
			$\bdry W = X \sqcup (-|X|)$.
		\end{enumerate}	
	\end{Lemma}
	
	\begin{Proof}
		Part (i): This is Lemma \ref{lemma:IntrinsicGluedCompletion}.
		Part (ii) is clear.
		
		Part (iii): 
		Define $W$ as the intrinsic gluing of $X \times [0,\onehalf]$ and $|X| \times (\onehalf,1]$
		(which has collared boundary at $\{1\}$ because $X$ has no strange points)
		along $X \times \{\onehalf\}$.
		
		Part (iv): For the $5$-gapped spaces themselves this is clear (their $4$-skeleton is a disjoint union of manifolds),
		for the intrinsic-gluings it follows from Lemma \ref{lemma:LocallyConelike}.
		
		Part (v): The stratification of intrinsic gluings was chosen so that the top-stratum is not \inquotes{interrupted}
		in the middle, so this is true by construction.
		
		By enforcing a dense top-stratum and disallowing codimension $1$-strata, orientations of the top-stratum of
		$|X|$ are in one-to-one correspondence with those of the top-stratum of $|X|$: By construction
		the top-stratum of the intrinsic gluing of $X \times [0,\onehalf]$ and $|X| \times (\onehalf,1]$
		is an open subset of the top-stratum of $|X| \times [0,1]$, and an orientation is obtained by restricting the orientation
		sheaf. For the other direction, note that an orientation of the intrinsic-gluing induces one of its boundary-component
		$|X| \times \{1\}$. If we reverse that orientation, we get an inverse construction to the one above.
	\end{Proof}
	
	\begin{Remark}
		For example when studying L-classes in §\ref{sec:LClasses}, we will be interested in
		bordism-invariants of subspaces obtained by cutting. If for example $X$ is a $5$-gapped MHSS with simple links etc.,
		we are interested in invariants of $X_{\{0\}}$ and of $(X_{\{0\}})_{\{0\}}$ and so on.
		Their transport under homeomorphism is then controlled by cutting the bordism $W$
		with $B_\pm$ the two-sides of $X_{\{0\}}$ (i.\,e.\ $B_+ = (X_\geq - \txt{boundary-collar}) \times I$)
		to get a cut of $|X|$ bordant to $X_{\{0\}}$ and so on.
		
		Note, that while the Lemma above looks like $W$ is essentially just a cylinder,
		this is \emph{not}
		true for this bordism of $X_{\{0\}}$:
		On the $|X|$-side boundary of $W$, the cut is obtained directly on $|X|$
		(so for example, $|X|_{\{0\}}$ never has strange-points, as this would destroy boundary-collars in $|X|_\geq$),
		so in general the cut $|X|_{\{0\}}$ is not simply
		the intrinsic space of $X_{\{0\}}$ and the bordism of Lemma \ref{lemma:BordismSpaceIntrinsic}
		fits the two together by a bordism, which has to have some \inquotes{twist} in the right-hand-side of $W$!
		
		In particular, this fitting together uses the relative version of Theorem \ref{thm:AlgTopCond}
		(our main source of examples) by using that the simple links hypothesis does not only
		kill end-obstructions, but also h-cobordism-obstructions, making cylinder-neighborhoods,
		thus cuts unique in a suitable sense.
	\end{Remark}
	
	\begin{Example}\label{example:Intrinsic5GappedSimpleLinks}
		(a) If $X$ has at most two meeting strata, then $|X|$ is $5$-gapped with simple links
		if $X$ is.
		
		This is, because that condition is local and only on pairs of meeting strata, and
		any pair of strata locally looks like a pair of strata in $X$ by (Q1).
		It is also clear from Rmk.\,\ref{rmk:IntrinsicConstruction}.
		
		(b) If $X$ is obtained from $|X|$ which is $5$-gapped with simple links --
		for example $|X|=M$ a manifold --
		by refining the stratification, e.\,g.\ to improve certain geometric properties (like making
		the space pl), and these changes do not break these properties on $X$, then 
		the results apply.
		This answers the question from the introduction about topological manifolds stratified as
		($5$-gapped with simple links) pl-pseudomanifolds: Indeed in such cases the subsequent results
		will apply, and Goresky--MacPherson L-classes and (poincaré-duals of) manifold L-classes coincide.
	\end{Example}
	
	\begin{Remark}
		In general, $X$ is $5$-gapped $\Rightarrow$ $|X|$ is $5$-gapped (again by (Q1) or Rmk.\,\ref{rmk:IntrinsicConstruction}).
		However, even if $X$ has simple links, this need not be true for $|X|$, where the issue occurs if
		middle strata are \inquotes{promoted away} (Rmk.\,\ref{rmk:IntrinsicConstruction}),
		for example if $| X \supset Y \supset Z | = (X \supset Z)$ it is not clear from the holink-fibers of $Z$ in
		$Y$ and in $X-Y$, what the holink-fibers of $Z$ in $X$ are.
		This is geometrically related to the problem of compositions of cylinder-neighborhoods not
		being cylinder-neighborhoods in general (see \citep[Example\,1 (p.\,5)]{EdwardsTopRegNbhs} for a counterexample).
		
		In this sense, the results as presented here relate the question about dependence on the
		choice of stratification with a question about local structure in MHSS.
		We cannot fully resolve this local structure question, but can demonstrate its accessibility
		e.\,g.\ for the case of spaces with only two strata.
		This is not an unusual restriction on results for the local structure of MHSS
		\citep{HTWW}.
		
		One might gain the impression, that a similar trick on the geometric side might work
		to show, that intrinsic stratifications of $5$-gapped MHSS with simple links automatically have a suitable cutting property:
		Requiring the existence of a $5$-gapped MHSS with simple-links $Y$
		such that $|Y| = X$, seems to work as a compatible condition on first sight
		(fixing a choice of $Y$ for the relative statements as for intrinsic gluings above).
		The lowest strata of $X=|Y|$ could be cut using manifold-transversality (so strange-points never occur),
		then extended using $Y$ (note, that the proof of Theorem \ref{thm:AlgTopCond} never actually uses
		that the boundary-collars of the cut of the lower skeleton are stratified) and \inquotes{simultaneously}
		by (Q3) to $X =|Y|$.
		However, the resulting $X_{\{0\}}$, while having no strange points, need not be a MHSS
		(we effectively replaced certain strata \inquotes{from below}, which may break holinks being fibrations).
		So cutting $X_{\{0\}}$ again may not work, as we know nothing about cylinder-neighborhood existence.
		
		An application where similar methods might conceivably work is products.
		While this paper is not primarily concerned with ring-structure on bordism-theories,
		this could be an interesting extension.
	\end{Remark}
	
	\begin{Remark}\label{rmk:CompatibleCategories}
		At this point, we want to briefly come back to the comment at the beginning of §\ref{sec:CompatibleConditions}
		concerning a more categorial formulation of compatibility:
		The results of the present section, mostly require the existence of suitable cutting properties,
		however,
		as we have seen in the last two sections, the relevant data required to give the transversality\Slash{}compatibility
		condition includes:
		\begin{itemize}
			\item The choice of (underlying) spaces of cycles.
			\item The careful choice of relative\Slash{}extendable cuts, e.\,g.\ boundary-com\-pat\-ible or
			middle-compatible.
			\item The careful choice of what extension means precisely (e.\,g.\ extending collars from $X^k$, but only
			spaces from $C$, see proof of Theorem \ref{thm:AlgTopCond}).
			\item The choice of relative regions (e.\,g.\ $X^k \cup C$ as compare to only the boundary in
			Lemma \ref{lemma:CutIntrinsicGluings}).
		\end{itemize}
		It seems, that at least all of this information needs to be captured by a systematic approach.
		While understanding the precise abstract structure behind these problems might lead to
		interesting insights, it seams unlikely that giving such a treatment here would considerably simplify things.
	\end{Remark}

	\section{Application to Singular L-Classes}\label{sec:LClasses}

	The construction of L-classes given by Thom \citep{ThomMexiko58}
	via signatures of suitable submanifolds has been studied and generalized extensively in the past.
	Most notably, it was a major influence that lead to the development of intersection homology theory
	\citep{GoreskyMacPherson} which provides bordism-invariant signatures for example for
	Siegel's Witt-spaces \citep{SiegelWittSpaces}.
	This Witt-condition and bordism-invariant signature were extended to MHSS by \citep{FriedmanPD},
	so are available for the theory of this paper.
	For previous results on the matter, see also §\ref{sec:DefsAndResults}.
	
	\begin{Definition}\label{def:Witt}
		By \citep[§8 (p.\,2197f)]{FriedmanPD}, there is a 
		\defName{Witt-condition}
		for MHSS (assuming compactly dominated local holinks as in \citep{HughesTeardrops,FriedmanPD},
		see Def.\,\ref{def:Basics}) with
		\inquotes{sufficiently many local approximate tubular neighborhoods} (we
		do not use any details of this local-neighborhood requirement, it is relatively weak and satisfied
		for the spaces used here, as detailed in the proof of Theorem \ref{thm:WittBordism} below).
		These MHSS Witt-spaces are endowed with a well-defined \defName{signature} of
		the (middle perversity, middle dimension) rational intersection-homology pairing.
	\end{Definition}
	
	The reference defines a Witt-condition for
	closed spaces, but the condition on links can be applied to spaces with boundary (because links in the boundary are
	the same as those in the collar) and the neighborhood requirement is not affected by this.
	What we refer to as Witt-space is an \inquotes{H$\Rationals$-Witt space} in 
	\citep{FriedmanPD}, i.\,e.\ a \inquotes{homotopy Witt-space} for the ring of rational numbers.
	
	\begin{Lemma}[Thom's Theorem for MHSS]\label{lemma:ThomForMHSS}
		This signature is Witt-bordism invariant, i.\,e.\ 
		given a compact Witt-space $X$, then
		\begin{equation*}			
			\sign(\bdry X) = 0
			\quad\txt{and}\quad
			\sign(X \sqcup (-Y)) = \sign(X) - \sign(Y)\txt{.}
		\end{equation*}
	\end{Lemma}
	\begin{Proof}
		The usual arguments (see for example \citep[p.\,155]{GoreskyMacPherson}
		or \citep[p.\,124f]{BanaglTopInvOfStratSp}) apply:
		Define $Y := X \cup_{\bdry X} \cone(\bdry X)$ by coning off the boundary, then	$Y$, away from the cone-point $v$,
		again has the property that inclusion of lower-middle perversity intersection-chains in
		upper-middle perversity intersection-chains induces an isomorphism on intersection-homology.
		(Technically this should be formulated on sheaves, where the correct statement is:
		The algebraic mapping cone of the inclusion-induced map lower-middle perversity inter\-sect\-ion-chain-sheaf in
		the upper-middle perversity intersection-chain-sheaf has non-trivial stalk only at $v$.)
		
		The cone-point also evidently has a locally conelike neighborhood, so the duality-results of
		\citep[Cor.\,7.3 (p.\,2193)]{FriedmanPD} apply, and the (homology) stalk at $v$ itself can be related to
		its link (which is $\bdry X$) by the cone-formula for intersection-homology.
		Combining both (duality and form of the homology at the link, which is the only non-vanishing contribution
		to the algebraic mapping-cone) one obtains a ladder-diagram as for the manifold-case
		(see e.\,g.\ \citep[p.\,124 \& p.\,127]{BanaglTopInvOfStratSp}),
		so the same algebraic arguments as for the proof of Thom's theorem on manifolds apply
		(the image of the map to he algebraic mapping-cone is a lagrangian-subspace).
	\end{Proof}
	
	\begin{Lemma}\label{lemma:IntrinsicGluingWitt}
		If $X$ is a MHSS Witt-space, then so is the intrinsic gluing (if defined, i.\,e.\ if $X$ has no strange points)
		$X \times [0,\onehalf] \GlueIntrin{X \times \{\onehalf\}} |X| \times (\onehalf,1]$.
	\end{Lemma}
	\begin{Proof}
		We defined Witt-spaces with boundary (Def.\,\ref{def:Witt}) to
		satisfy the same local holink condition that \citep{FriedmanPD} uses for closed spaces.
		In particular $X \times S^1$ (which has the same local holinks as $X$) is H$\Rationals$-Witt in the
		sense of \citep{FriedmanPD}.
		This condition can be phrased (see \citep[p.\,2198]{FriedmanPD}) independently of the choice of stratification
		via topological invariance of the used intersection-homology condition \citep{QuinnIntrinsicSkeleta}.
		This is seen by the uniqueness of the Deligne-sheaf satisfying certain axioms (see e.\,g.\ \citep{FriedmanPD}).
		These axioms fix, as part of the definition an orientation-sheaf (or local-coefficient system)
		on the top-stratum, so one should be careful, not to introduce a non-trivial choice here,
		for example one may not introduce a codimension $1$ stratum.
		This is not an issue here (Lemma \ref{lemma:BordismSpaceIntrinsic}\,v, see also its proof)
		by construction of the intrinsic gluing.
		In particular, gluing $X \times [0,\onehalf] \GlueIntrin{X \times \{\onehalf\}} |X| \times (\onehalf,1]$
		to $-(X \times [0,\onehalf] \GlueIntrin{X \times \{\onehalf\}} |X| \times (\onehalf,1])$
		(i.\,e.\ to itself with orientation reversed),
		which is a stratification of the underlying space of $X \times S^1$, is Witt.
		Again, since this is a local condition via the local holinks, the intrinsic gluing
		$X \times [0,\onehalf] \GlueIntrin{X \times \{\onehalf\}} |X| \times (\onehalf,1]$
		itself is Witt.
	\end{Proof}
	
	Using this, we can prove Theorem \ref{thm:WittBordism} about the existence of an oriented bordism-theory
	of certain (\inquotes{admissible}, see below) MHSS-Witt-spaces.
	
	\begin{ProofOfThm2}
		In this proof, we call spaces \inquotes{admissible}, if they are intrinsic multi-gluings of $5$-gapped MHSS
		with simple links (Def\,\ref{def:MultiGluing}), that further have dense top-stratum
		and satisfy the Witt-condition.
		MHSS Witt-spaces never have codimension $1$ strata, see \citep[p.\,2197]{FriedmanPD}, so
		admissible spaces have no codimension $1$ stratum.
		We will call a space $X$ \inquotes{super-admissible}, if $|X|$ is admissible and $X$ has no strange points.
		
		We first check, that being admissible is a compatible condition:
		By Lemma \ref{lemma:IntrinsicGluedCompletion}, this is the case for being an intrinsic multi-gluings of $5$-gapped
		MHSS with simple links.
		Having dense top-stratum and no codimension $1$ stratum is auxiliary by Lemma \ref{lemma:ExamplesOfCompatibleConditions}
		thus can be added to obtain another compatible condition by Lemma \ref{lemma:IntercompatibilityOfConditions}.
		For the local holink condition in the definition of Witt-spaces (Def.\,\ref{def:Witt}),
		the same argument applies.
		
		Witt-spaces are additionally required to satisfy a local neighborhood condition.
		This is satisfied automatically -- thus the compatible condition given above is equivalent to 
		being admissible -- as follows:
		By Lemma \ref{lemma:BordismSpaceIntrinsic}\,iv, all spaces involved have locally conelike $4$-skeleton.
		Local conelike neighborhoods are mapping-cylinder-neighborhoods thus teardrop-neighborhoods or	
		\inquotes{local approximate tubular neighborhoods}, so $X^4$ always has 
		\inquotes{sufficiently many local approximate tubular neighborhoods}.
		By Hughes' neighborhood theorem \citep{HughesNbhEx02},
		$X$ itself has \inquotes{sufficiently many local approximate tubular neighborhoods} if $X^4$ does,
		i.\,e.\ in high dimensions, their existence is automatic.
		
		Thus being admissible is a compatible condition and we can apply
		Theorem \ref{thm:Bordism} to obtain an oriented (the addition of orientation does not pose a problem
		either by Lemma \ref{lemma:BordismSpaceIntrinsic}\,v) bordism-theory $\Omega^{\txt{MHSS-Witt}}_*$
		of admissible spaces,
		which is a generalized homology-theory and
		provides stratified orientation-preserving stratified homeomorphism invariant fundamental classes for closed oriented
		admissible spaces,
		in particular for closed oriented $5$-gapped MHSS Witt-spaces with simple links and dense top-stratum
		(Lemma \ref{lemma:BordismSpaceIntrinsic}\,ii).
		
		Given a (general\Slash{}unstratified) homeomorphism of super-admissible spac\-es $h: X \rightarrow Y$,
		by Lemmas \ref{lemma:BordismSpaceIntrinsic} and \ref{lemma:IntrinsicGluingWitt},
		the intrinsic gluing $X \times [0,\onehalf] \GlueIntrin{X \times \{\onehalf\}} |X| \times (\onehalf,1]$
		is an admissible bordism (over the point) $X \bordant |X|$.
		This clearly has underlying topological space $X \times I$, so is, via projection to the $X$ coordinate,
		an admissible bordism over the underlying space of $X$ as well, so $[X] = [|X|] \in \Omega^{\txt{MHSS-Witt}}_n(X)$,
		with $n = \dim(X)=\dim(Y)$.
		Analogously one obtains $[Y] = [|Y|] \in \Omega^{\txt{MHSS-WittWitt}}_n(Y)$.
		By property (Q1) of Quinn's intrinsic stratifications (see §\ref{sec:IntrinsicTransversality}),
		the homeomorphism of the underlying spaces $h: X \rightarrow Y$ is stratified with respect to the intrinsic stratifications
		$h : |X| \rightarrow |Y|$ thus by Theorem \ref{thm:Bordism}\,ii $h_*([|X|]) = [|Y|]$.
		All put together, thus $h_*([X]) = h_*([|X|]) = [|Y|] = [Y]$.
		
		By Lemma \ref{lemma:ThomForMHSS} the signature (of cycles) induces a group-homo\-mor\-phism
		$\sigma : \Bordism_*(A) \rightarrow \Integers, [f: X \rightarrow A] \mapsto \sign(X)$ (see
		Theorem \ref{thm:Bordism}\,iii), which,
		by definition on fundamental-classes $[X] := [\id : X \rightarrow X]$ is normalized s.\,t.\ $\sigma([X]) = \sign(X)$.
	\end{ProofOfThm2}
	
	We show, that given this bordism-theory of Witt-spaces, one can define, for spaces $X$ with a fundamental class $[X]$ in
	that bordism-theory, certain
	mappings $l_{n-4k} : \pi^{n-4k}(X) \rightarrow \Integers$
	(for $4 k < (n-1)/2$, when the cohomotopy-sets are actually groups, see the stabilization-remark in the proof of
	Lemma \ref{lemma:ltoL}), which using the construction of \citep{GoreskyMacPherson}
	induce homological L-classes (see Lemma \ref{lemma:ltoL}).
	
	Usually (e.\,g.\ in \citep{GoreskyMacPherson}) these $l_{n-4k}$
	are constructed by making a representative $\varphi : X \rightarrow S^{n-4k}$ transverse to
	a base-point $* \in S^{n-4k}$ and then mapping $[\varphi] \mapsto \sign( \varphi_\perp^{-1}(\{*\}) )$.
	This clearly requires a transversality result, and products with intervals (for well-definedness).
	These properties are implicit with the bordism-theory that was constructed above:
	In particular, repeated application of inverse of the excision-isomorphism (or similarly,
	the desuspension homomorphism of the associated reduced homology-theory)
	models this kind of transversality with trivial normal-bundle.
	Hence we define:
	
	\begin{Definition}
		Define $l_k : \pi^{n-4k}(X) \rightarrow \Integers$ as the mapping
		that assigns to $[\varphi] \in \pi^{n-4k}(X)$,
		the integer $\sigma( \desusp^{n-4k}( \varphi_*([X]) ) )$,
		with $\desusp$ denoting the inverse of the suspension-homomorphism of the reduced theory
		$\tilde\Omega^\mathcal{C}_*$ associated to $\Omega^\mathcal{C}_*$.
		
		This is simply the signature $\sign(Y) = \sigma([Y])$ of the geometric cycle
		(the domain $Y$ of a representative $[f:Y \rightarrow X]$) of
		$[f] := \desusp^{n-4k}( \varphi_*([X]) ) \in \tilde\Omega^\mathcal{C}_{4k}(S^0)
		= \Omega^\mathcal{C}_{4k}(\{\txt{pt}\})$.
	\end{Definition}
	
	The desuspension map can be constructed either from, or analogously to, the inverse of the excision-isomorphism
	(see §\ref{sec:Bordism}).
	
	\begin{Remark}\label{rmk:SmallLInvariant}
		Since these depend only on the fundamental class $[X]$,
		their trans\-port-behavior (e.\,g.\ under homeomorphisms)
		is determined by the trans\-port-behavior of those fundamental-classes.
	\end{Remark}
	
	\begin{Lemma}\label{lemma:ltoL}
		These $l_*$ (when defined) induce homological L-classes as in \citep{GoreskyMacPherson}.
		If additionally $X$ is a pl-pseudo\-mani\-fold, then these L-classes are the Goresky--MacPherson L-classes.
		
		These L-classes are invariant under stratified orientation-preserving ho\-meo\-morphisms,
		and under orientation-preserving general (unstratified) ho\-meo\-morphisms of spaces with no strange points,
		if the intrinsic stratification of the underlying space is $5$-gapped with simple links.
	\end{Lemma}
	
	\begin{Proof}
		The construction of homological L-classes
		(via potentially stabilizing with spheres to get to high degrees, tensoring with the rationals,
		using a theorem of Serre that rational co-homotopy stably becomes co-homology,
		and the universal-coefficient-theorem to get from
		dual cohomology to homology) works exactly as in \citep{GoreskyMacPherson}
		(see also \citep[§5.7 (p.\,120--122)]{BanaglTopInvOfStratSp}).
		
		If $X$ is a pl-pseudomanifold, it can be cut as a pl-pseudomanifold using the pl-structure.
		This cut is in particular also topologically a cut (i.\,e.\ a cut by Def.\,\ref{def:Cut}). By the construction of desuspension,
		either analogously to the proof of excision\Slash{}Lemma \ref{lemma:Excision}, or from excision, it is
		clear, that well-definedness (see proof of Lemma \ref{lemma:Excision}) implies, that this pl-cut is an allowed choice
		also in the topological theory.
		Desuspending $n-4k$ times, thus cutting $n-4k$ times, is the same as making a map to a $S^{n-4k}$ transverse to a point,
		and taking signature of the component over the base-point of $S^0$
		thus is the same as the construction of \citep{GoreskyMacPherson},
		which is taking the signature of a transverse-preimage.
		
		The invariance properties follow by Rmk.\,\ref{rmk:SmallLInvariant}.
	\end{Proof}

	\appendix

\section{Properties of Intrinsic Spaces and Gluings}

\renewcommand{\thesection}{A}

The intrinsic gluings used in §\ref{sec:IntrinsicTransversality} are essentially
half-intrinsic cylinders (see e.\,g.\ \citep{FriedmanBordismPseudomanif} for a similar construction
in the pl-case). This appendix shows that
they are indeed MHSS and preserve local-conelikeness of $4$-skeleta.
The expected enumeration of all relevant cases of meeting strata and
corresponding elementary arguments works, details are included for the sake of completeness, nothing
unexpected happens.

\begin{Lemma}\label{lemma:LocallyConelike}
	Local conelikeness:
	\begin{enumerate}[(i)]
		\item If $X^k$ is locally conelike, and $Y$ is a coarsening of $X$, then $Y^k$ is locally conelike.
		In particular $|X|^k$ is then locally conelike.
		\item If $X^k$, $Y^k$ and $Z^{k-1}$ are locally conelike, $k\geq 4$, then $(X \GlueIntrin{Z} Y)^k$ is locally conelike.
		The intrinsic gluing $X \GlueIntrin{Z} Y$ is defined in Def.\,\ref{def:CoarseGluing}.
	\end{enumerate}
\end{Lemma}

\begin{Proof}
	Part (i):
	Let $y \in Y_i$, then $y \in X_j$ for $j \leq i < k$ and there are a compact $L$,
	an open neighborhood $U$ of $y$ in $X^k$ and
	$\varphi : \opencone{L} \times \Reals^j \rightarrow U$ mapping $(v,0) \mapsto y$.
	As stratified homeomorphism, $\varphi$ induces a one-to-one mapping
	\begin{align*}
		\phi: \{ \txt{components of strata of $L$} \} &\xrightarrow{\txt{1:1}} \{ \txt{components of strata of $U-X_j$} \}	\\
		S\qquad &\makebox[\widthof{$\xrightarrow{\txt{1:1}}$}][c]{$\mapsto$}
		\qquad \varphi\big( S \times (0,1) \times \Reals^j \big)
	\end{align*}
	and by $Y$ being a coarsening of $X$, also $U'$ (induced by $Y$ on the underlying space of $U$)
	is a coarsening of $U$ and thus induces a mapping
	\begin{equation*}		
		\psi : \{ \txt{components of strata of $U$} \} \rightarrow \{ \txt{components of strata of $U'$} \}\txt{.}
	\end{equation*}
	Stratify the underlying space of $L$ as $L'$ with strata the pre-images of strata of $U'$
	under $\psi \circ \phi$.
	Then $\varphi : \opencone{L'} \times \Reals^j \rightarrow U'$,
	is again a stratified homeomorphism by construction.
	
	Part (ii):
	This local condition need be checked near $Z^{k-1}$ only, thus let $z \in Z_j$, where $j \leq k-1$.
	If $j = k-1$
	
	If $j < k-1$, we simplify notation by writing $Z=\opencone(L) \times \Reals^j$ with elements $z=([l,s], p)$,
	for $U \subset Z^{k-1}$, and as (unstratified) topological space $W = Z \times (-2,2)$
	for a neighborhood of $z$ in $(X \GlueIntrin{Z} Y)^k$ (near the gluing-interface the underlying space has
	this product-form by definition).
	In this local standard-form, the construction can be achieved in the \inquotes{plane} spanned by
	$t \in (-2,2)$ and the cone-coordinate $t$ of $Z$ as follows:
	
	Let $q$ the quotient map $q : L \times [0,1) \times \Reals^j \times (-2,2) \rightarrow W$ of the cone,
	and define $L_W$ as follows:
	First, let $D$ the \inquotes{diamond}-shape
	$D := \{ (l,s,t) \in L \times [0,1) \times (-2,2)$~$|$ $s=|t| \}$ and
	the underlying space $L_W := q(D)$ of what will become the link.
	Since $k\geq 4$, by (Q2) $|Z \times \Reals|$ is the topologically intrinsic stratification,
	thus homogeneous along the the $\Reals$-coordinate, so we may
	pick any $\tau \in \Reals$ and the subspace $Z' := Z \times \{\tau\} \subset |Z \times \Reals|$
	does not depend on this choice $\tau$, and $|Z \times \Reals| = Z' \times \Reals$ up to stratified
	homeomorphism. $Z \times \{\tau\} \neq |Z|$ in general, and Z' may not be a MHSS near zero-strata of $|Z|$.
	For $0\leq m<k-j-1$ (recall: We want to show the $k$-skeleton is conelike, and are constructing the
	conelike-neighborhood for a point $z$ in the $j$-stratum of $Z$, so $k-j-1$ is the expected formal link-dimension,
	in particular $L_W^{k-j-1} := L_W$),
	assume at first $k<n$, where $n = \dim(X)=\dim(Y)$ (the top-stratum is not in play), under the convention
	that $L^a = \emptyset$ for $a < 0$:
	\begin{align*}
		L_W^m = \{v\} \times \Gamma
		&\cup q( L^m \times \{ \onehalf \} \times \{0\} )\\
		&\cup q( D \cap L^{m-1} \times (0,1) \times (-2,2) )
	\end{align*}
	where either $\Gamma = \{ \pm 1 \}$ if $z \in Z'_j$ (if the point of interest is in the same stratum of the coarsening
	$Z'$ as in $Z$) or $\Gamma = \{-1 \}$ if $z\in Z'_{j'}$ for a $j'>j$.
	If $k = n$, the top-stratum becomes relevant, and we change the above definition for $k-j-2$
	(the highest skeleton not yet equal to the full $L_W$) to
	\begin{equation*}
		L_W^{k-j-2} = \{v\} \times \Gamma \cup q( D \cap L^{k-j-3} \times (0,1) \times (-2,2) )\txt{.}
	\end{equation*}
	Finally we need the stratified homeomorphism to a neighborhood of $z$, thus define
	(where $[l,s] \in \opencone(L)$, $t \in (-2,2)$ with $([l,s],t) \in q(D)$,
	further $s'$ is the cone-coordinate of $\opencone(L_W)$,
	and $p\in \Reals^j$; the result is written as an element of $\opencone(L) \times (-2,2) \times \Reals^j$)
	\begin{align*}
		\varphi_W : \opencone(L_W) \times \Reals^j &\rightarrow W	\\
		([[l,s],t,s'],p)			&\mapsto ( [l,s s'], t s',p )\txt{.}
	\end{align*}
	On underlying spaces, this is just the product of $\Reals^j$ with the cone on $D$ times $L$, so
	$\varphi_W$ is certainly a homeomorphism to its image.
	Note, that the cone-point $v'$, where $s'=0$, is mapped to $(v, 0, p)$ and we identified coordinates
	using $\varphi \times \id_{(-2,2)}$, which maps $(v, p=0) \mapsto z$, so $\varphi_W( v', 0 ) = z$ as required.
	Further $\varphi_W$ essentially interpolates by the cone-coordinate $s'$ from the \inquotes{diamond} $D$
	inwards to $(v, 0, 0)$, which is stratified, because along $s=0$ on the $t<0$ side this is along the $j$-stratum of
	$X$, on the $t>0$ side it is either also along the $j$-stratum (if present in the coarsening $Z'$ and
	$\Gamma = \{ \pm 1 \}$) or this is part of a higher-stratum of $Y$, but then also $+1 \not\in \Gamma$
	so this is also reflected by the definition of the link.
	Along $t=0$ this is just the cone-neighborhood of $z\in Z$. In between we stay in a single stratum of
	$Z-Z_j\times (-2,0)$ or $Z-Z_j\times (0,2)$ when moving along the cone-coordinate in the range $(0,1)$,
	where no change of stratum occurs either, by definition of the stratification on the cone.
\end{Proof}

\begin{Lemma}\label{lemma:HalfIntrinsicCyl}
	Given a compact MHSS $Z^n$, then $W$ with underlying space
	$W^{n+1} = Z \times (0,1)$ filtered for $k<n+1$ by
	\begin{equation*}
		W^k := Z^{k-1} \times (0,\onehalf] \cup |Z \times (\onehalf,1)|^k \cup
		\begin{cases}
			\emptyset & \txt{if }k\geq n	\\
			Z^k \times \{\onehalf\}	& \txt{otherwise}
		\end{cases}
	\end{equation*}
	is a MHSS.
\end{Lemma}

\begin{Proof}
	First, note that strata are manifolds: For $k<n+1$ they are disjoint unions of
	strata of the form $|Z \times (\onehalf,1)|_k$, $Z_{k-1} \times (0,\onehalf)$ or $Z_k \times \{\onehalf\}$,
	which are products of manifolds, because $Z$ and $|Z|$ have manifold-strata.
	The top-stratum is an open subset of $|Z|_n \times (0,1)$, which is a manifold, thus it is also a manifold.
	Further, the filtration is evidently by closed subsets and certainly $W$ is metric and separated because $X$ is.
	
	To show, that $W$ is a MHSS, since this is a local condition, and the left-hand-side (lhs) of $W$ is $X\times(0,\onehalf)$
	and the rhs is $|X|\times (\onehalf,1)$ which are MHSS,
	it suffices to show that $W$ is a MHSS near the middle at $\onehalf$.
	Pairs of components of strata, $S, S'$, that meet here, are (up to reordering) of the form $S$ is a component of
	$X_n \times \{\onehalf\}$
	and $\dim(S') > \dim(S)$.
	We distinguish four cases:
	
	(i) Case $S' \subset X \times \{\onehalf\}$:
	Since $X \times \{\onehalf\}$ is a MHSS, this is tame and has the holink a fibration.
	
	(ii) Case $S' \subset \txt{lhs}$: Clearly $S' = S \times (0,\onehalf)$,
	thus we have $(S \cup S', S) = (S \times (0,\onehalf], S\times \{\onehalf\})$, which is
	certainly tame and has the holink a fibration.
	
	(iii) Case $S' \subset \txt{rhs}$:
	If $S' = S \times (\onehalf,1)$, this is essentially the same as case (ii).
	Otherwise, we first note, that $S' \subset |X \times (\onehalf, 1)|$, which is a coarsening of $X \times (\onehalf, 1)$
	by (Q1).
	So $S'$ has a stratification $Y$ as subspace of $X \times (\onehalf, 1)$ and hence as MHSS.
	By definition of the product-stratification,
	$S \times (\onehalf,1) \subset Y$ is a stratum. We have a partial ordering on strata of $Y$ by
	$Y_\alpha \geq Y_\beta :\Leftrightarrow \closure(Y_\alpha) \supset Y_\beta$ (because the stratification of $Y$ is
	induced from a filtration by closed subsets). Extend this partial order to a total order to get
	$Y_\alpha \geq Y_\beta \geq \ldots \geq Y_\omega = S \times (\onehalf,1) \geq \ldots$.
	Since $Y$ is a MHSS, there is a tameness retraction $R_\alpha$ of $Y_\alpha \cup Y_\beta$ to $Y_\beta$ and so on.
	Now $R := R_\alpha * R_\beta * \ldots$ is well-defined, because $R_\alpha(x,1) \in Y_\beta$ (the domain of the next
	term $R_\beta$),
	and $R(x,1) \in Y_\omega = S \times (\onehalf,1)$.
	This, in turn, has again a tameness-retraction as in case (ii).
	The composition is nearly-strict, because $Y_\alpha, \ldots, Y_\omega \subset S'$.
	The holink is a fibration, because $R$ induces a fibered deformation-retraction of the holink of $(S \cup S', S)$
	to $(S \times [\onehalf,1), S \times \{\onehalf\})$.
	
	(iv) Case $S'$ is (a component of) the top-stratum of $W$: This is the only remaining case,
	as the top-stratum is the only stratum that is not in lhs, middle or rhs.
	Note that on the (closed) lhs of the top-stratum $W_{n+1}$ is simply $X_n \times (0,\onehalf]$.
	We want to push $S'$ into the (open) rhs and then apply the argument of the previous case (iii)
	as follows:
	First, define $T' : X_n \times (0,1) \times I \rightarrow X_n \times (0,1)$
	as either $T'(x,s,t) := (x, t s + (1-t) \onehalf)$ if $s \leq \onehalf$ or
	$T'(x,s,t) := (x, s)$ if $s \geq \onehalf$. This is clearly continuous, can be extended
	by the identity on the rhs and moves things in the (closed) rhs $\pi_2 T'(x,s,1) \in [\onehalf, 1)$.
	Further, to push things into the open rhs $|X| \times (\onehalf, 1)$,
	while still having a continuous extension by the identity, use the metric distance
	$d : X_n \times (0,1) \rightarrow (0,1], y \mapsto \min(1, \dist(y, X^{n-1} \times [\onehalf, 1]))$,
	where the distance is measured in the metric space $X\times [0,1]$.
	This is well-defined, because the space $X^{n-1} \times [\onehalf, 1]$ is compact.
	Let $S(s,d) := d  \frac{3+s}{4}	+ (1-d) s$, which is chosen such that for $s \in [\onehalf,1)$ and
	$d \in (0,1]$ both $S(s, d \searrow 0) \rightarrow s$ and $S(s,d) > \onehalf$.
	Hence $T'' : X_n \times [\onehalf,1) \times I \rightarrow X_n \times [\onehalf,1)$ given by
	$T''((x,s),t) := (x, (1-t)s + t S(s, d(x,s))$ is such that $\pi_2 T''((x,s),1) \in (\onehalf, 1)$
	is in the (open) rhs,
	while for $d \searrow 0 : T''((x,s),t) \rightarrow (x,s)$, so that $T''$ can be extended by the identity on
	the remainder of $X \times [\onehalf, 1)$.
	If we define $T := T' * T''$ extended by the identity to $W_{n+1} \cup S$,
	then $X_n \times (0,1)$ gets pulled into $X_n \times (\onehalf, 1)$ in a stratified way,
	because $|X\times (\onehalf,1)|$ being a coarsening of $X\times (\onehalf,1)$ implies $X_n \times (0,1) \subset W_{n+1}$.
	Now, the argument of case (iii) applies.
\end{Proof}

\begin{Lemma}\label{lemma:CoarseGluingMHSS}
	The intrinsic gluing of Def.\,\ref{def:CoarseGluing} $X \GlueIntrin{Z} Y$ is a MHSS.
\end{Lemma}

\begin{Proof}
	Since the boundaries are collared by definition, and being a MHSS is a local condition,
	the claim needs to be checked only on open boundary-collars, where the stratification is the one
	of the Lemma \ref{lemma:HalfIntrinsicCyl} given above, thus a MHSS.
\end{Proof}

\section*{Acknowledgements}

The results presented here were obtained throughout the author's phd-thesis at Heidelberg University
carried out under the supervision of Markus Banagl,
to whom I am grateful for many suggestions and much support provided towards this research.\\
The author was personally supported by the \inquotes{Studienstiftung des Deutschen Volkes} through
a scholarship.\\
I further want to thank the anonymous reviewer for a very thorough reading of the manuscript and
very helpful comments.

	\bibliographystyle{elsarticle-num} 
	\bibliography{refs.bib}
	
	
		
		
		
\end{document}